\documentclass[12pt]{amsart}

\usepackage{graphicx}
\usepackage{amscd}

\usepackage{epsfig}
\usepackage{times}
\usepackage{psfig}

\RequirePackage{ifthen}
\RequirePackage{times}
\RequirePackage{mathptm}
\RequirePackage{pifont}

\newtheorem{theorem}{Theorem}
\theoremstyle{plain}

\newtheorem{lemma}{Lemma}

\newtheorem{proposition}{Proposition}

\numberwithin{equation}{section}

\newcommand\C{{\mathbb C}}
\newcommand\R{{\mathbb R}}

\begin{document}
\title[Euler polynomials]{On the zero attractor of the Euler polynomials}
\author{William M.Y.  Goh}
\address{Department of Mathematics, Drexel university, Philadelphia, PA}
\email{wgoh@math.drexel.edu}
\author{Robert  Boyer}
\address{Department of Mathematics, Drexel university, Philadelphia, PA}
\email{rboyer@math.drexel.edu}
\date{August 31, 2004}

\begin{abstract}
We study the limiting behavior of the zeros of the Euler polynomials.
When linearly scaled, they approach a definite curve in the complex plane
related to the Szeg\"o curve 
which governs the behavior of the roots of
the Taylor polynomials associated to the exponential function.
Further, under a conformal transformation, the scaled zeros are uniformly distributed.
\end{abstract}

\maketitle

\section{Introduction}

Eighty years ago Szeg\"{o}  \cite{szego} studied the asymptotic behavior of
the roots of the Maclurin polynomials associated with the exponential
function. He found that if the roots are linearly scaled relative
to the degree then the roots
approach a curve $S$ (see Figure 2)  in the complex plane
given by  $z \in \C$ such that $|z e^{1-z}|=1$ and 
$|z| \leq 1$. The behavior of the roots and 
poles of the Pad\'{e} approximants and other Taylor polynomials have been
analyzed \cite{wielonsky}.

On the other hand,
given any sequence of polynomials $\{p_{n}(x)\}$, where $p_{n}(x)$ is of
degree $n$,  asking how are the zeros of $p_{n}(x)$ distributed
in the complex $x$-plane is too  general to get a reasonable
answer. The best we can hope for is to focus on a special family of
polynomials where a definite answer is possible.

In this paper, we initiate the study of the asymptotic behavior of the roots
of  Euler polynomials
$E_{n}(x)$ which are defined by means of generating functions as
\begin{equation}
\frac{2e^{\xi x}}{e^{\xi }+1}=\sum_{n\geq 0}E_{n}(x)\frac{\xi ^{n}}{n!}.
\label{euler-po}
\end{equation}
This generating function is listed
among the principal generating functions by Louis Comtet 
\cite{comtet} for combinatorial applications.
Their linearly scaled roots approach a curve related to the Szeg\"{o}
curve $S$
together with an  interval on the real axis.

Polynomials of binomial
type were introduced by Rota and Mullin  \cite{ro&mu}. The reason for the name
is that if 
\[
e^{yD(x)}=\sum_{n\geq 0}\frac{\phi _{n}(y)}{n!}x^{n},
\]
where $D(x)$ is a polynomial, it follows that

\[
\phi _{n}(u+v)=\sum_{r}(_{r}^{n})\phi _{r}(u)\phi _{n-r}(v),
\]
which is a reminiscent of the binomial theorem. 
Herbert Wilf has the same
description in his book
``generatingfunctionology''
\cite{wilf}. 
Strictly speaking, $E_{n}(x)$ is not a polynomial of binomial
type. However, it is close to being binomial type and has the simplest $D(x)$
function. In this case $D(\xi )=\xi$.
We hope to investigate the zero attractors of this wider class of polynomial
families in the future.  As evidence, the family of Bernoulli polynomials, for example,
are easily handled with the techniques in this paper.
A study of the behavior of their real zeros was recently done \cite{veselov}.

\medskip
Let $\{p_{n}(x)\}$ be a sequence of polynomials.  A set $A$ in the $x$-plane
is called the 
\textit{zero attractor of zeros}  of $\{p_{n}(x)\}$ if the following two
conditions hold:

{
\leftskip=12pt

\noindent
a) Let $A_{\varepsilon }:=\bigcup _{x\in A} \, B(x,\varepsilon )$, where 
$B(x,\varepsilon )$ is the open disc centered at $x$ with radius $\varepsilon 
$. That is, \ $A_{\varepsilon }$ is just the $\varepsilon $-neighborhood of
the set $A$. Then there exists an integer $n_{0}(\varepsilon )$, for all 
$n\geq n_{0}$, all zeros of $p_{n}(x)$ are in $A_{\varepsilon }$.

}

{
\leftskip=12pt

\noindent
b) For all $x\in A$ and for all $\varepsilon >0$, there exists an integer 
$n_{1}(x,\varepsilon )$ and a zero $r$ of the polynomial $p_{n_{1}}(x)$ such
that $r\in B(x,\varepsilon )$.

}  

\medskip

\noindent 
Condition b) simply says that every point of $A$ is an
accumulation point of zeros of $\{p_{n}(x)\}$.

The Euler polynomials $E_{n}(x)$ are defined in  (\ref{euler-po}).  Since the
nearest singularity to the origin of $\frac{1}{e^{\xi }+1}$ are $\xi =\pm
\pi i,$ it is easy to see that for all $x\in {\mathbb {C}},$ the power series in
(\ref{euler-po})  converges absolutely and uniformly on any compact subset in 
$\left| \xi \right| <\pi$. In other words, although the polynomial $E_{n}(x)$
is defined for all complex $x$ but the power series is convergent only for
 $\xi $ with $\left| \xi \right| <\pi$.  By the Cauchy residue theorem, we
have:

\[
\frac{E_{n}(x)}{n!}=\frac{2}{2\pi i}
\oint_{\left| \xi \right| =1}
\frac{e^{x\xi }}{(e^{\xi }+1)\xi ^{n+1}} \, d\xi
\]
This integral expression is valid for all $x\in \C$. Let $x$ be replaced by 
$nx$ and we can write the above equation as:
\[
\frac{E_{n}(nx)}{n!}
=\frac{2}{2\pi i}
\oint_{\left| \xi \right| =1}
\left( \frac{e^{x\xi }}{\xi }  \right)^{n}
\, \frac{ 1 }{\xi (e^{\xi }+1)} \, d\xi
\]
The goal of this paper is to study the zero distribution of the polynomial 
$E_{n}(nx)$.

\section{A Generalization of the Szeg\"{o} Approximation}

We state two generalizations of the Szeg\"{o} approximation. Let 
\[
S_{n}(z) :=  \sum_{j=0}^{n}\frac{z^{j}}{j!}.
\]

\begin{proposition}
\label{szego1}Let $S$ be a subset contained in $\left| z\right| >1$ so that
the distance between $S$ and the unit circumference $\left| z\right| =1$ is 
$\delta >0$ and let $\alpha $ be chosen so $1/3<\alpha <1/2$.  Then 
\begin{equation*}
\frac{S_{n-1}(nz)}{e^{nz}}=
\frac{(ze^{1-z})^{n}}{\sqrt{2\pi n}(z-1)}
\left(1+O(n^{1-3\alpha }) \right),
\end{equation*}
where the constant in the big $O$ term is uniform for all $z\in S$.
\end{proposition}

\begin{proof}
By residue theory, we have for $R>0$ 
\begin{equation*}
S_{n-1}(z)=
\frac{1}{2\pi i}
\oint_{\left| \zeta \right| =R}
\frac{e^{\zeta }}
{\zeta -z}\frac{\zeta ^{n}-z^{n}}   {\zeta ^{n}}  \, d\zeta,
\end{equation*}
Note that $\zeta =z$ is a removable singularity of the integrand. Therefore,
the above expression is valid for all complex $z$.  For asymptotics in the
region $\left| z\right| \geq 1+\delta$,  we choose the contour to be the
circle $\left| \zeta \right| =n$. Thus 
\[
S_{n-1}(nz)=\frac{1}{2\pi i}
\oint_{\left| \zeta \right| =n}
\frac{e^{\zeta }}
{\zeta -nz}\frac{\zeta ^{n}-(nz)^{n}}{\zeta ^{n}}  
\, d\zeta ,
\]
Since $nz$ is not included inside the contour $\left| \zeta \right| =n$,  a
simple application of Cauchy's Theorem gives: 

\begin{eqnarray*}
S_{n-1}(nz)
&=&
\frac{-(nz)^{n}}{2\pi i}
\oint_{\left| \zeta \right| =n}
\frac{e^{\zeta }}{(\zeta -nz)\zeta ^{n}} \, d\zeta  \\
&=&
\frac{-z^{n}}{2\pi i}\oint_{\left| \zeta \right| =1}
\frac{e^{n\zeta }}
{(\zeta -z)\zeta ^{n}}  \, d\zeta   \\
&=&
\frac{-z^{n}}{2\pi i}
\oint_{\left| \zeta \right| =1}
\frac{e^{n(\zeta -\ln \zeta )}}{\zeta -z}  \, d\zeta,
\end{eqnarray*}

\noindent
where $\ln \zeta $ is the principal branch with $-\pi <\arg \zeta \leq \pi $
and $n$ is a positive integer. Now we apply the saddle point method to
construct the asymptotics of the integral. Since the critical point is 
$\zeta =1,$ the neighborhood of $\zeta =1$ must be carefully analyzed. Let 
$\eta =n^{-\alpha }$ with $1/3<\alpha <1/2$.  The contour integral is
decomposed into two integrals: 
\[
\frac{1}{2\pi i}
\oint_{\left| \zeta \right| =1}
\frac{e^{n(\zeta -\ln \zeta )}}{\zeta -z}  \,  d\zeta =I_{1}+I_{2},
\]
where $I_{1}$ is the integral on the circular arc in a small neighborhood of 
$1:  \,  -\eta \leq \arg \zeta \leq \eta ,\left| \zeta \right| =1$,  and the
integral $I_{2}$ is along the path in the remaining part of the circle.   Set 
$\zeta =e^{i\theta }$ in $I_{1}$.  The Taylor expansion of the integrand in a
small neighborhood of $\theta =0$ is worked out below: 
\[
e^{n(\zeta -\ln \zeta )}=e^{n}e^{n(-\frac{\theta ^{2}}{2}+
O\left(\theta^{3}\right))}
=e^{n}e^{-n\theta ^{2}/2} \left(1+O(n^{1-3\alpha }) \right).
\]
Inserting these estimates in $I_{1}$ and carrying out some simplifications
we get 
\begin{eqnarray*}
I_{1}&=&\frac{1}{2\pi i}\int_{-\eta }^{\eta }
\frac{e^{n}e^{-n\theta^{2}/2}(1+O(n^{1-3\alpha }))}
{e^{i\theta }-z}ie^{i\theta } \, d\theta  \\
&=&
\frac{e^{n}}{2\pi }\int_{-\eta }^{\eta }
\frac{e^{-n\theta^{2}/2}(1+O(n^{1-3\alpha }))}{1-z+O(n^{-\alpha })}
\left(1+O(n^{-\alpha }) \right)  \, d\theta
\end{eqnarray*}
Since $1/3<\dot{\alpha}<1/2,$ we have $3\alpha -1<\alpha$.
 So the error
term $O(n^{-\alpha })$ is absorbed into $O(n^{1-3\alpha })$.  Hence,  for 
$\left| z\right| \geq 1+\delta$,  we see 
\[
I_{1}=\frac{e^{n}}{2\pi (1-z)} 
\left(\int_{-\eta }^{\eta }
e^{-n\theta^{2}/2}  \, d\theta \right)
\, \left(1+O(n^{1-3\alpha }) \right),
\]
where the big $O$ term holds uniformly for $\left| z\right| \geq 1+\delta$.
If we put $n\theta ^{2}/2=u^{2}$,  we get 
\[
I_{1}=\left(\frac{e^{n}}{2\pi (1-z)}\sqrt{\frac{2}{n}}
\int_{-\omega }^{\omega}  e^{-u^{2}} \, du \right)  \, \left(1+O(n^{1-3\alpha })\right),
\]
where $\omega =\sqrt{ \frac{n^{1-2\alpha }}{2}  }$.
 Since $\alpha <1/2$, $\omega$
tends to $\infty $ with $n$.  But, as 
$n\rightarrow \infty$, 
$\displaystyle
\int_{\omega}^{\infty }e^{-u^{2}} \,  du
=O(e^{-\omega ^{2}}/\omega )=o(n^{1-3\alpha })$.  We
may therefore replace the limits of integration by $\pm \infty $ without
altering the error term $1+O(n^{1-3\alpha })$. This gives 
\[
I_{1}= \left(\frac{e^{n}}{2\pi (1-z)}\sqrt{\frac{2}{n}}
\int_{-\infty }^{\infty}  e^{-u^{2}}  \, du\right)\,  \left (1+O(n^{1-3\alpha })\right) 
=
\frac{e^{n}}{\sqrt{2\pi n}(1-z)}(1+O(n^{1-3\alpha })),
\]
after taking $\int_{-\infty }^{\infty }  e^{-u^{2}} \, du=\sqrt{\pi }$ into
consideration. To justify that $I_{1}$ gives the major contribution we
obtain an upper estimate for $I_{2}$:
\begin{equation*}
\left| I_{2}\right| \leq \frac{1}{2\pi }\int_{ C }
\frac{e^{n\Re(\zeta )}}{  \left| \zeta -z\right|   }
 \, \left| d\zeta \right| ,
\end{equation*}
where $C$ is the contour determined by $\eta \leq \left| \arg \zeta \right|
\leq \pi ,$ and $\left| \zeta \right| =1$.  Obviously, 
$\Re(\zeta)\leq \cos \eta $. Note that $\left| \zeta -z\right| \geq \delta $ for
all $z\in S$.  Hence 
\begin{equation*}
\left| I_{2}\right| \leq \frac{1}{2\pi }
\frac{e^{n\cos \eta }(2\pi )}{\delta }=
\frac{e^{n\cos \eta }}{\delta }.
\end{equation*}
Upon using 
\[
\cos \eta =1-\frac{\eta ^{2}}{2}+O(\eta ^{4}),
\]
we see that 
$e^{n\cos \eta }=e^{n}e^{-\frac{1}{2}n^{1-2\alpha}}
(1+O(n^{1-3\alpha }))$.  But the factor 
$e^{-\frac{1}{2}n^{1-2\alpha }}
=o(\frac{1}{\sqrt{n}})$.  Consequently, 
\begin{equation*}
I_{2}=o(I_{1}).
\end{equation*}
This completes the proof of Proposition \ref{szego1} .
\end{proof}

Next, the following proposition states the asymptotics of $S_{n}(z)$ in the
region $\Re(z)<1$.

\begin{proposition}
\label{szego2}
For $1/3<\alpha <1/2$,  we have 
\begin{equation*}
\frac{S_{n-1}(nz)}{e^{nz}}=1-\frac{(ze^{1-z})^{n}}{\sqrt{2\pi n}(1-z)}
\left(1+O(n^{1-3\alpha }) \right),
\end{equation*}
where the big $O$ constant holds uniformly for an arbitrary compact set 
$K\subseteq \Re(z)<1$.
\end{proposition}

Thus the ordinary Szeg\"{o} approximation is generalized from the open disc 
$\left| z\right| <1$ to the open half plane $\Re(z)<1$. 

\begin{proof}
The proof is actually very similar to that of Proposition \ref{szego1}. Here
we use a suitable integral representation for $S_{n-1}(z)$ in the region. 
We start off with: 
\begin{equation*}
S_{n-1}(z)=\frac{1}{2\pi i}\oint_{C}\frac{e^{\zeta }}{\zeta -z}
\frac{\zeta^{n}-z^{n}}{\zeta ^{n}} \, d\zeta,
\end{equation*}
where $C$ is any closed contour encircling the origin.  This integral
representation is valid for all complex $z$. We insert $nz$ for $z$  to
obtain 
\begin{equation*}
S_{n-1}(nz)=\frac{1}{2\pi i}\oint_{C}\frac{e^{\zeta }}{\zeta -nz}
\frac{\zeta^{n}-(nz)^{n}}{\zeta ^{n}} \, d\zeta.
\end{equation*}
Let $K$ be an arbitrary compact set in the domain $\Re(z)<1$.  We
choose the circular integration contour 
$C:\zeta =-R+(R+1) e^{i\theta }, \,  -\pi \leq \theta \leq \pi $ 
with $R$ so large that $K$ becomes strictly included
inside $C$. \ Let $z$ belong to $K$ and replace the contour $C$ by $nC$.
This gives 
\begin{equation*}
S_{n-1}(nz)=\frac{1}{2\pi i}
\oint_{nC}\frac{e^{\zeta }}{\zeta -nz}
\frac{ \zeta ^{n}-(nz)^{n}}{\zeta ^{n}} \, d\zeta .
\end{equation*}
After a change of variables $\zeta \rightarrow n\zeta $ and the use of the
Cauchy Integral Formula, we get 
\begin{equation*}
S_{n-1}(nz)=e^{nz}-\frac{(nz)^{n}}{2\pi i}
\oint_{C}
\frac{e^{n\zeta }}{\zeta -z}\frac{d\zeta }{(n\zeta )^{n}}.
\end{equation*}
This implies 
\begin{equation*}
\frac{S_{n-1}(nz)}{e^{nz}}=1-\frac{e^{-nz}z^{n}}{2\pi i}
\oint_{C}
\frac{e^{n(\zeta -\ln \zeta )}}{\zeta -z}  \, d\zeta.
\end{equation*}
where $C$ is the circle $\zeta =-R+(R+1)e^{i\theta },
\, -\pi \leq \theta \leq \pi$.  
The critical point is
still $\zeta =1$ which is a point in the contour $C$.   It is clear that the
saddle point method can be applied to the contour integral.  The procedure
is very similar to what we did in the previous case. We omit the details.
This is how we prove the statement of Proposition \ref{szego2}. 
\end{proof}

We can rederive the same result by using the integral representation of 
$S_{n}(z)$ in the original derivation of the Szeg\"{o} approximation. 
Szeg\"{o} \cite{varga}  used 
\begin{equation*}
\frac{S_{n}(nz)}{e^{nz}}=1-\frac{n^{n+1}e^{-n}}{n!}
\int_{0}^{z}(ve^{1-v})^{n} \, dv
\end{equation*}
to obtain the result of the proposition in the disc $\left| z\right| <1$. 
A good elaboration of his approach can indeed offer a proof of Proposition
\ref{szego2}. The saddle point method as we presented in the proof is just
an easier way to get the result.

\begin{proposition}
\label{szego3}(Jet Wimp) \ (Uniform\ Approximation of $S_{n}(nt)$ for 
$t\geq 0)$
\begin{equation*}
\frac{S_{n}(nt)}{e^{nt}}
=\delta (t)+\sqrt{\frac{2}{\pi }}\frac{\xi (t)t}{t-1}
{\rm Erfc}(\sqrt{n}\xi (t))
\,   \left(1+O(\frac{1}{\sqrt{n}}) \right)
\end{equation*}
uniformly for $t\geq $ $0$, where $\xi (t)=\left| t-1-\ln t\right| ^{1/2}$,
and 
\begin{equation*}
\delta (t)=\left\{ 
\begin{array}{rcl}
1,&  \mathrm{for}& 0 \leq t  <  1, \\ 
0,&  \mathrm{for}& t  \geq 1.
\end{array}
\right.
\end{equation*}
and for all $t, {\rm Erfc}(t):=\int_{t}^{\infty }e^{-s^{2}} \, ds$.
\end{proposition}

This version of the Szeg\"{o} approximation was proved by Jet Wimp (personal
communication).

\section{A Decomposition for the Euler Polynomial}

The Euler polynomial $E_{n}(nx)/n!$ is decomposed as the sum of two
polynomials:

\begin{proposition}
\label{m-kprop}Let $\mu $ be an integer $\geq $ 0 and let $F_{\mu }(\xi )$
be given as 
\begin{equation*}
F_{\mu }(\xi ):=\frac{1}{\xi (e^{\xi }+1)}+\sum_{k=0}^{\mu }
\left[\frac{1}{(2k+1)\pi i(\xi -(2k+1)\pi i)}+
\frac{1}{-(2k+1)\pi i(\xi +(2k+1)\pi i)} \right].
\end{equation*}
\noindent
Then we have 
\begin{equation*}
\frac{E_{n}(nx)}{n!}=M_{n,\mu }(x)+K_{n,\mu }(x),
\end{equation*}
where 
\begin{eqnarray*}
M_{n,\mu }(x)
&=& \frac{1}{\pi i}\oint_{\left| \xi \right| =1}
\left(\frac{e^{x\xi }}{\xi } \right)^{n}F_{\mu }(\xi ) \, d\xi,
\\
K_{n,\mu }(x)
&=& 2\sum_{k=0}^{\mu }
\left[\frac{S_{n-1}(nx(2k+1)\pi i)}
{((2k+1)\pi  i)^{n+1}}+\frac{S_{n-1}(-nx(2k+1)\pi i)}{(-(2k+1)\pi i)^{n+1}} \right].
\end{eqnarray*}
Here $S_{n}(z):=\sum_{j=0}^{n}z^{j}/j!$,  the $n$-th partial sum of $e^{z}$
as usual.
\end{proposition}

\begin{proof}
The following integral representation for $E_{n}(x)$ is valid for all 
$x\in \C$: 
\begin{equation*}
\frac{E_{n}(x)}{n!}=\frac{2}{2\pi i}\oint_{\left| \xi \right| =1}
\frac{ e^{x\xi }}{(e^{\xi }+1)\xi ^{n+1}}  \, d\xi .
\end{equation*}
Let $x$ be replaced by $nx$ to get 
\begin{equation*}
\frac{E_{n}(nx)}{n!}=\frac{2}{2\pi i}
\oint_{\left| \xi \right| =1}
\left(\frac{e^{x\xi }}{\xi }\right)^{n}
\frac{  d\xi }{\xi (e^{\xi }+1)}.
\end{equation*}
Note that for each integer $k\geq 0$, 
\begin{equation*}
-\frac{1}{(2k+1)\pi i(\xi -(2k+1)\pi i)}  -  
\frac{1}{(-(2k+1)\pi i)(\xi +(2k+1)\pi i)}
\end{equation*}
is the sum of singular parts of $\frac{1}{\xi (e^{\xi }+1)}$ at $(2k+1)\pi i$
and $-(2k+1)\pi i$. Hence $F_{\mu }(\xi )$ is analytic in the annulus 
$0<\left| \xi \right| <(2\mu +3) \pi $. 
The bigger $\mu $ is, the bigger the
domain of analyticity of $F_{\mu }(\xi )$ is.  Hence 
\begin{eqnarray}
\lefteqn{
\frac{E_{n}(nx)}{n!}
=
\frac{2}{2\pi i}
\oint_{\left| \xi \right| =1}
\left(\frac{e^{x\xi }}{\xi }\right)^{n}   F_{\mu }(\xi )  \,  d\xi   
- 
\frac{2}{2\pi i}  \oint_{\left| \xi \right| =1}
\left(\frac{e^{x\xi }}{\xi } \right)^{n}        
}  \nonumber \\
&&        
\quad \times
\left(\sum_{k=0}^{\mu }
\left[\frac{1}{(2k+1)\pi i(\xi -(2k+1)\pi i)} + 
\frac{1}{(-(2k+1)\pi i)(\xi +(2k+1)\pi i} \right] \right) 
\, d\xi  \label{m-kequ}
\end{eqnarray}

A typical term of the above sum is: 
\begin{eqnarray*}
\lefteqn{
\frac{2}{2\pi i}\oint_{\left| \xi \right| =1}\left( \frac{e^{x\xi }}{\xi }  \right)^{n}
\frac{1}{(2k+1)\pi i(\xi -(2k+1)\pi i)}  \, d\xi  }         \\
&&
=
\frac{2}{2\pi i}\frac{1}{(2k+1)\pi i}\oint_{\left| \xi \right| =1}
\frac{-(\frac{e^{x\xi }}{\xi })^{n}}{(2k+1)\pi i}
\sum_{l=0}^{\infty }(\frac{\xi }{ (2k+1)\pi i})^{l}  \,  d\xi 
\end{eqnarray*}
after using a geometric series expansion. Obviously, the series is uniformly
convergent on $\left| \zeta \right| =1$,  we carry out the integration term
by term. 
By the Cauchy integral theorem only those terms with $l\leq n-1$
survive. Thus we get 
\begin{eqnarray*}
\lefteqn{
\frac{1}{\pi i}  \frac{(-1)  } {  ((2k+1)\pi i)^2   }
\sum_{l=0}^{n-1}
\
\left( \frac{1}{  (2k+1)\pi i} \right)^{l}    
\,   \oint_{  \left| \xi \right| =1  } \, e^{nx\xi }  \xi ^{l-n} \, d\xi  
}              
\\
&&
\quad=
\frac{1}{\pi i}
\frac{(-1)}{((2k+1)\pi i)^{2}}   \sum_{l=0}^{n-1}
(\frac{1}{(2k+1)\pi i})^{l}\frac{(nx)^{n-l-1}}{(n-l-1)!}(2\pi i)
\\
&&
\quad=
 \frac{(-2)}{((2k+1)\pi i)^{2}}\sum_{l=0}^{n-1}
\frac{(nx(2k+1)\pi i)^{n-l-1}}{(n-l-1)!}((2k+1)\pi i)^{-n+1} \\
&&
\quad=
(-2)((2k+1)\pi i)^{-n-1}
\sum_{j=0}^{n-1}\frac{(nx(2k+1)\pi i)^{j}}{j!}  \\
&&
\quad=
\frac{(-2)}{((2k+1)\pi i)^{n+1}        }
             S_{n-1}(nx(2k+1)\pi i),
\end{eqnarray*}
where $S_{n}(z)$ is the $n$-th partial sum of $e^{z}$; that is, 
$S_{n}(z)=\sum_{j=0}^{n}z^{j}/j!$. \ In a similar way, we obtain 
\begin{eqnarray*}
\lefteqn{
\frac{2}{2\pi i}\oint_{\left| \xi \right| =1}
 \left( \frac{e^{x\xi }}{\xi }  \right)^{n}
\frac{1}{(-(2k+1)\pi i)(\xi +(2k+1)\pi i)} \, d\xi 
}         \\
&&
\quad \quad=
\frac{-2}{(-(2k+1)\pi i)^{n+1}}S_{n-1}(-nx(2k+1)\pi i).
\end{eqnarray*}
Inserting these back into  (\ref{m-kequ}), we complete the proof of the
proposition.
\end{proof}

\noindent 
 Proposition \ref{m-kprop}  is important since it provides
asymptotics for the Euler polynomials in various regions.  The asymptotics
for $M_{n,\mu }(x)$ (for any fixed $\mu$)  can be easily found by the
classical saddle point method.  It is 
\[
M_{n,\mu }(x)=\sqrt{\frac{2}{\pi }}
\left((xe)^{n}F_{\mu }(1/x)\frac{1}{\sqrt{n}x} \right) 
\, \left(1+O(\frac{1}{n}) \right),
\]
  where the big 
$O$ term holds uniformly for 
$\varepsilon \leq \left| 1/x\right| \leq (2\mu +1)\pi +\varepsilon$.
The uniformity can be justified using the fact that 
$F_{\mu }(\xi )$ is analytic in $0<\left| \xi \right| <(2\mu +3)\pi$.  The
asymptotics for $K_{n,\mu }(x)$ can be obtained from the generalized 
Szeg\"{o}'s approximation.  Using these asymptotic approximations we will prove 
that the point set $K$ defined below is the zero attractor of
the Euler polynomials and can also  determine their density distribution.

\section{The Zero Attractor}

Let the point set $K$ be defined by the graph in Figure 1:

\begin{figure}[ht]
\centerline{\psfig{figure=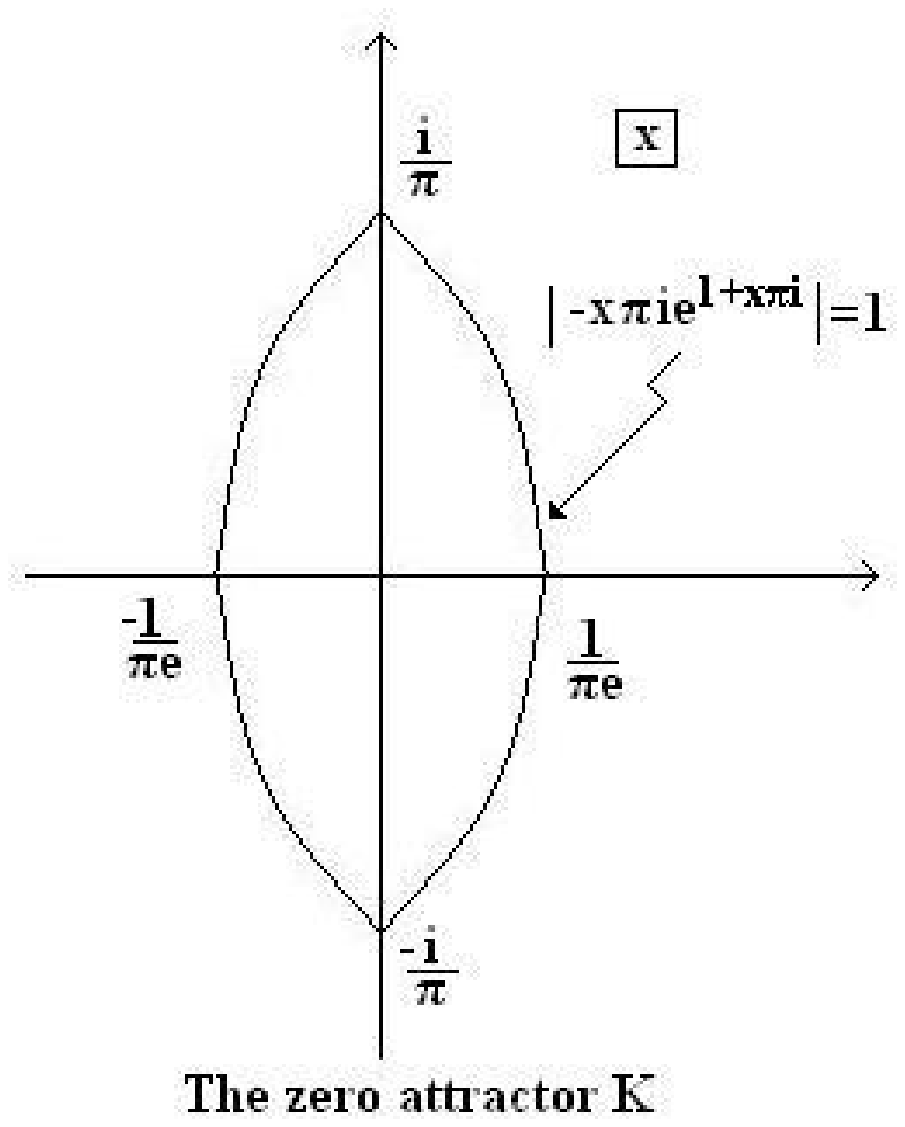, height=3.5in, width=4.0in}}
\caption{}
\label{fig:zero_attractor_1}
\end{figure}

The point set $K$ consists of the curves as indicated and the real interval 
$[-\frac{1}{\pi e},\frac{1}{\pi e}]$.  $K$ is symmetric with respect to the
imaginary axis in the $x$-plane.  We will show as the first step that $K$
contains all accumulation points of the zeros of $E_{n}(nx)$.  In the
second step we will carry out the density calculations of the zeros residing
in an immediate neighborhood of $K$.  As a consequence, this will establish
that every point $K$ is an accumulation point of zeros of $E_{n}(nx)$.
Hence $K$ is precisely the zero attractor of $E_{n}(nx)$.  The density
calculation gives also the statistical distribution of zeros which includes
the information about the fraction of zeros along each segment of $K$.  The
work of proving that $K$ contains all accumulation points is divided into
the three parts  in the following lemma:

\begin{lemma}
\label{zeroattra}
(a) Let $\left| x_{0}\right| >\frac{1}{\pi },$ then $x_{0}$
is not an accumulation point of zeros of \ the Euler polynomials $E_{n}(nx)$.

\noindent
(b) There is no accumulation point of zeros in the region 
$\left| x\right|  \leq \frac{1}{3\pi }$ except for real numbers.

\noindent
(c) If $x$ is a non-real accumulation point of zeros in the region 
$\frac{1}{3\pi }<\left| x\right| \leq \frac{1}{\pi },$ then we must have either 
\begin{equation*}
\left| x\pi ie^{1-x\pi i}\right| =1
\quad
\textrm{or} 
\quad
\left| -x\pi ie^{1+x\pi i}\right| =1.
\end{equation*}
\end{lemma}

\begin{proof}
Let us prove part (a). Now 
\begin{eqnarray*}
\frac{E_{n}(nx)}{n!}
&=&
\frac{2x^{n}}{2\pi i}\oint_{\left| \xi \right| =1}
\left(\frac{e^{\xi }}{\xi } \right)^{n}\frac{d\xi }{\xi (e^{\xi /x}+1)}  \\
&=&
\frac{2x^{n}}{2\pi i}\int_{\left| \xi \right| =1}\exp [nf(\xi )] 
\frac{d\xi }{\xi (e^{\xi /x}+1)},
\end{eqnarray*}
where $f(\xi )=\xi -\ln \xi $. We choose the principal branch of $\ln \xi $
here. To invoke the saddle point method, we need to find the critical points
which are roots of $f^{\prime }(\xi )=0$ or $1-\frac{1}{\xi }=0$, that is, 
$\xi =1$ is the only critical point in question. Observe that 
$\left| \exp f(\xi )\right| =\left| \frac{e^{\xi }}{\xi }\right| =e^{\cos \theta }$,
where $\xi =e^{i\theta }$. Now $e^{\cos \theta }$ attains its maximum at 
$\theta =0,$ hence at $\xi =1$. By the saddle point method  \cite{cop} we have 
\begin{eqnarray*}
\frac{E_{n}(nx)}{n!}
&=& \frac{2x^{n}}{2\pi i}  
\left(    \frac{e^{n}}{e^{1/x}+1}
\left(\frac{-2\pi }{nf^{\prime \prime }(1)}\right)^{1/2} \,
\left(1+O(\frac{1}{n}) \right) \right)  \\
&=&
\frac{2(ex)^{n}}{\sqrt{2\pi n}}\frac{1}{e^{1/x}+1}
\left(1+O(\frac{1}{n}) \right),
\end{eqnarray*}
where the big $O$ constant can be made uniform for a given compact set 
$K\subseteq \{x:\left| x\right| >\frac{1}{\pi }\}$. From this, we see that
for any given set $K\subseteq \{x:\left| x\right| >\frac{1}{\pi }\}$, there
exists $n_{0}$ such that for all $n\geq n_{0}$ the polynomial $E_{n}(nx)$
has no zeros in $K$. Hence the proof of part (a). 

Next we prove part (b). Assume that there is a non-real accumulation point
of zeros in the region $\left| x\right| \leq \frac{1}{3\pi },$ say $x_{0}$.
 Of course, $x_{0}\neq 0$, then there exists an integer $\mu _{0}\geq 2$
and a positive number $\varepsilon >0$ such that 
\begin{equation*}
\frac{1}{(2\mu _{0}+1)\pi }+\varepsilon <\left| x_{0}\right| \leq 
\frac{1}{(2\mu _{0}-1)\pi }
\end{equation*}
Since $x_{0}$ is an accumulation point of zeros of $E_{n}(nx)$, there exists
an infinite sequence of integers $n_{j}$ such that 
\begin{equation*}
E_{n_{j}}(n_{j}x_{n_{j}})=0\text{ and }x_{n_{j}}\rightarrow x_{0}
\text{ as }  j\rightarrow \infty .
\end{equation*}
We may assume that for all large $j$,  $\{x_{n_{j}}\}$ is in the region 
\begin{equation*}
\frac{1}{(2\mu _{0}+1)\pi }+\varepsilon <\left| x_{0}\right| \leq 
\frac{1}{ (2\mu _{0}-1)\pi }+\varepsilon ,
\end{equation*}
so we apply Proposition  \ref{m-kprop} with $\mu $ chosen as $\mu _{0}-1$ and
keep $\mu _{0}$ fixed in the following arguments. Use the asymptotics of 
$M_{n,\mu }(x)$ we see that $E_{n_{j}}(n_{j}x_{n_{j}})=0$ implies 
\begin{equation}
\sqrt{\frac{2}{\pi }}\frac{(x_{n_{j}}e)^{n_{j}}}{\sqrt{n_{j}}x_{n_{j}}}
F_{\mu _{0}-1}(\frac{1}{x_{n_{j}}}) \left(1+O(\frac{1}{n_{j}}) \right)+
K_{n,\mu_{0}-1}(x_{n_{j}})=0,  \label{m-kinzero}
\end{equation}
where the big $O$ term holds uniformly for all $x$ in the region 
\begin{equation*}
\frac{1}{(2\mu _{0}+1)\pi }+\varepsilon <\left| x\right| \leq 
\frac{1}{(2\mu_{0}-1)\pi }+\varepsilon .
\end{equation*}
Note that the summation in $K_{n,\mu _{0}-1}(x_{n_{j}})$ is running from $0$
to $\mu _{0}-1$ and for all large $n_{j}$%
\begin{eqnarray*}
\left(  \frac{1}{(2\mu _{0}+1)\pi }+\varepsilon \right) (2k+1)\pi
 &<&
\left|  x_{n_{j}}(2k+1)\pi i\right| \\
&\leq& 
\left(\frac{1}{(2\mu _{0}-1)\pi }+\varepsilon \right)  (2k+1)  \pi 
\end{eqnarray*}

A typical term in $K_{n,\mu _{0}-1}(x_{n_{j}})$ is 
\begin{equation*}
\frac{S_{n_{j}-1}(n_{j}x_{n_{j}}(2k+1)\pi i)}{((2k+1)\pi i)^{n_{j}+1}}.
\end{equation*}
Introducing $t=x_{n_{j}}(2k+1)\pi i$. For $0\leq k\leq \mu _{0}-2$,  we have 
$\left| t\right| \leq 1-\varepsilon$.
 Hence we apply Szeg\"{o} approximation
of Proposition  \ref{szego2}. Thus 
\begin{eqnarray*}
\lefteqn
{
\frac{S_{n_{j}-1}(n_{j}t)}{((2k+1)\pi i)^{n_{j}+1}}=
\frac{e^{n_{j}t}}{ ((2k+1)\pi i)^{n_{j}+1}}
} \\
&& \qquad
\times \, \left[ 1-\frac{1}{\sqrt{2\pi n_{j}}}\frac{1}{1-t}
(te^{1-t})^{n_{j}} \left(1+O(n_{j}^{1-3\alpha }) \right)  \right].
\end{eqnarray*}

The big $O$ constant in the above approximation is uniform since 
$\left| t\right| \leq \frac{2\mu _{0}-3}{2\mu _{0}-1}+\varepsilon$.
 Introduce the
function 
\[
g(x):=\frac{e^{\pi xi}}{e\pi xi}, \, \textrm{ or }  =(\pi xie^{1-\pi xi})^{-1}.
\]
Thus the above equation when expressed in term of $g(x)$ becomes 
\begin{eqnarray}
\lefteqn{
\frac{S_{n_{j}-1}(n_{j}t)}{((2k+1)\pi i)^{n_{j}+1}}=
\frac{(ex_{n_{j}})^{n_{j}}g^{n_{j}}((2k+1)x_{n_{j}})}{(2k+1)\pi i}
\cdot 
}  \nonumber \\
&&
\quad \times \,
\left[
 1-\frac{1}{\sqrt{2\pi n_{j}}}\frac{1}{1-t}
g^{-n_{j}}((2k+1)x_{n_{j}}) \left(1+O(n_{j}^{1-3\alpha }) \right)
\right],  \label{k-term}
\end{eqnarray}
where $1/3<\alpha <1/2$. Since $x_{0}$ is not real, so we may assume that 
$x_{0}$ is in the lower half plane. Thus $x_{0}=r_{0}e^{-i\theta _{0}}$ with 
$0<\theta _{0}<\pi $ and $\frac{1}{(2\mu _{0}+1)\pi }+\varepsilon <r_{0}
\leq 
\frac{1}{(2\mu _{0}-1)\pi }$.  Now, 
\begin{equation*}
\left| g(x_{0})\right| =\frac{e^{\pi r_{0}\sin \theta _{0}}}{e\pi r_{0}}
>\left| g(-x_{0})\right| =\frac{e^{-\pi r_{0}\sin \theta _{0}}}{e\pi r_{0}}.
\end{equation*}
Furthermore, $\frac{e^{\pi r\sin \theta _{0}}}{e\pi r},$ as a function of 
$r$ attains its minimum at $r=\frac{1}{\pi \sin \theta _{0}}$.  For each 
$1\leq k\leq \mu _{0}-2$ (remember $\mu _{0}\geq 2)$ we have 
\begin{equation*}
r_{0}<(2k+1)r_{0}<(2\mu _{0}-1)r_{0}\leq \frac{1}{\pi }\leq 
\frac{1}{\pi \sin \theta _{0}}.
\end{equation*}
Hence 
\begin{equation}
\left| g(x_{0})\right| >\left| g((2k+1)x_{0})\right| \geq \left|
g(-(2k+1)x_{0})\right| .  \label{g-inq}
\end{equation}
Also, since $x_{0}$ lies in the interior of the Szeg\"{o}'s domain 
\begin{equation*}
\left| z_{0}e^{1-z_{0}}\right| \leq 1\text{ and }\left| z_{0}\right| \leq 1.
\end{equation*}
Here, $z_{0}=\pi x_{0}i$. So, 
\begin{equation*}
\left| g(x_{0})\right| >1.
\end{equation*}
This prepares the behavior of the terms in $K_{n,\mu _{0}-1}(x_{n_{j}})$
with $0\leq k\leq \mu _{0}-2$. The term corresponding to $k=\mu _{0}-1$ in
the sum $K_{n,\mu _{0}-1}(x_{n_{j}})$ is 
$\frac{S_{n_{j}-1}(n_{j}x_{n_{j}}
(2\mu _{0}-1)\pi i)}{((2\mu _{0}-1)\pi i)^{n_{j}+1}}$ which can be estimated
using Proposition  \ref{szego3}. Let 
$t_{\mu _{0}}=\frac{n_{j}}{n_{j}-1} x_{n_{j}}(2\mu _{0}-1)\pi i$ 
and apply the proposition to get 
\begin{eqnarray}
\lefteqn{
\left| S_{n_{j}-1}(n_{j}x_{n_{j}}(2\mu _{0}-1)\pi i)\right| 
=\left|
S_{n_{j}-1}((n_{j}-1)t_{\mu _{0}})\right| 
}  \nonumber \\
&&
\quad
\leq S_{n_{j}-1}((n_{j}-1)\left| t_{\mu _{0}}\right| )
\nonumber \\
&&
\quad
\leq e^{(n_{j}-1)  \left| t_{\mu _{0}}\right| }
\,  2\sqrt{\frac{2}{\pi }}
\frac{\xi (\left| t_{\mu _{0}}\right| )\left| t_{\mu _{0}}\right| }
{\left| t_{\mu_{0}}\right| -1}
\,
{\rm Erfc}\left(\sqrt{n_{j}-1}\xi (\left| t_{\mu _{0}} \right| ) \right)
\left(1+O(\frac{1}{\sqrt{n_{j}}})\right) . \label{szego3-1}
\end{eqnarray}

\noindent
Note that $\left| t_{\mu _{0}}\right| \leq (\frac{1}{(2\mu _{0}-1)\pi }
+\varepsilon )(2\mu _{0}-1)\pi =1+(2_{\mu _{0}}-1)\pi \varepsilon $ as 
$n_{j}\rightarrow \infty$.  Recall that 
\begin{equation*}
g(x_{0})=e^{\pi x_{0}i}/e\pi x_{0}i,
\end{equation*}
where 
\begin{equation*}
x_{0}=\left| x_{0}\right| e^{-i\theta o}, \, 0<\theta _{0}<\pi
\end{equation*}
and 
\begin{equation*}
\frac{1}{(2\mu _{0}+1)\pi }+\varepsilon <\left| x_{0}\right| 
\leq \frac{1}{(2\mu _{0}-1)\pi },  \,\,  \mu _{0}\geq 2.
\end{equation*}
A detailed study of the Szeg\"{o}'s curve, defined as 
$\left| ze^{1-z}\right| =1$ and $\left| z\right| \leq 1$,  shows the following
features (see Figure 2):
\begin{figure}[ht]
\centerline{\psfig{figure=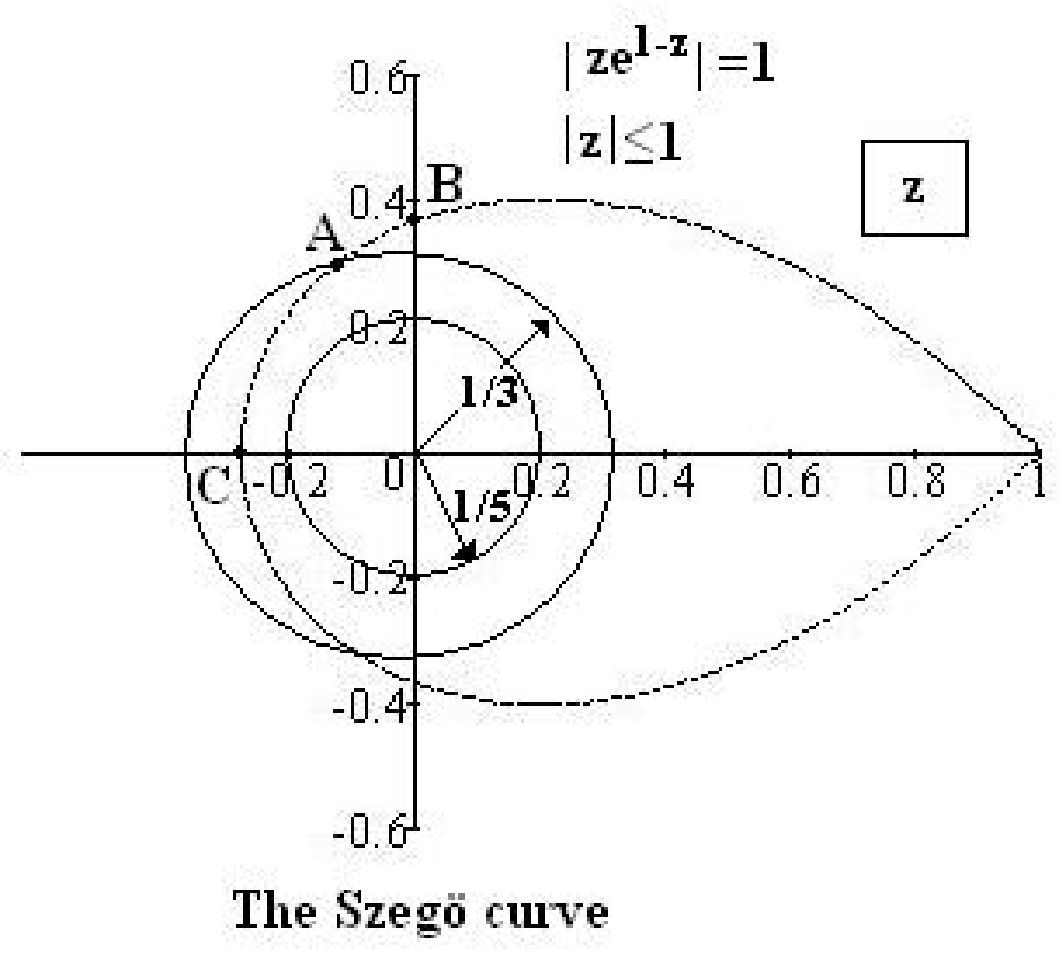, height=3.5in, width=4.0in}}
\caption{}
\end{figure}

\noindent
point $B$ is the intersection point of the curve with the imaginary axis and 
equals
$i/e,$
point $C$ is the intersection point of the curve with the negative real axis  
and is
$\approx -0.278\cdots$   (the unique real root of $e^{x-1}=-x)$
point $A$ is  the intersection point of
curve with the circle $\left| z\right| =1/3$. \ $\Re(A)  \approx
-9.861\times 10^{-2}$ that comes from solving the real root of 
$1/9-x^{2}=e^{2x-2}-x^{2}$ for $x$.
Hence, $z_{0}=\pi x_{0}i$ falls in the
interior of the Szeg\"{o}'s curve (also to the right of the imaginary axis
in the $z$-plane).  This implies that 
\begin{equation*}
\left| g(x_{0})\right| =\frac{1}{z_{0}e^{1-z_{0}}}>1.
\end{equation*}
Recall $\xi (\left| t_{\mu _{0}}\right| )=
\left| \left| t_{\mu _{0}}\right|
-1-\ln \right| t_{\mu _{0}}||^{1/2}$. \ So it is easy to see that 
\begin{equation*}
\frac{\xi (\left| t_{\mu _{0}}\right| )\left| t_{\mu _{0}}\right| }{\left|
t_{\mu _{0}}\right| -1}=O(1)
\end{equation*}
and  
$\displaystyle
{\rm Erfc}(\sqrt{n_{j}-1}\xi (\left| t_{\mu _{0}}\right| )):=
\int_{\sqrt{n_{j}-1}\xi (\left| t_{\mu _{0}}\right| )}^{\infty }e^{-s^{2}}ds=O(1)$
uniformly.
Hence from  (\ref{szego3-1})  we get 
\begin{equation}
\left| S_{n_{j}-1}(n_{j}x_{n_{j}}(2\mu _{0}-1)\pi i)\right| \leq
Ke^{n_{j}(1+\varepsilon )},  \label{s-esti}
\end{equation}
for some absolute constant $K$.  Now we are ready to see a contradiction
from  (\ref{m-kinzero}). Dividing (\ref{m-kinzero})  by $(x_{n_{j}}e)^{n_{j}}$
we get 
\begin{equation}
\sqrt{\frac{2}{\pi }}\frac{F_{\mu _{0}-1}
(1/x_{n_{j}})}{\sqrt{n_{j}}x_{n_{j}}}
\left(1+O(\frac{1}{n_{j}})\right)+(ex_{n_{j}})^{-n_{j}}K_{n,\mu _{0}-1}(x_{n_{j}})=0
\label{m-kinzero1}
\end{equation}
By  (\ref{s-esti})  the term with summation index $k=\mu _{0}-1$ in 
$(ex_{n_{j}})^{-n_{j}}K_{n,\mu _{0}-1}(x_{n_{j}})$ is estimated as 
\begin{eqnarray}
\left| \frac{S_{n_{j}-1}(n_{j}x_{n_{j}}(2\mu _{0}-1)\pi i)}{(ex_{n_{j}}
(2\mu_{0}-1)\pi i)^{n_{j}}}\right| 
&\leq&    K \left(\frac{e^{1+\varepsilon }}
{e\left| x_{n_{j}}\right| (2\mu _{0}-1)\pi } \right)^{n_{j}} \nonumber  \\
&=&
K\left(\frac{e^{\varepsilon }}{\left| x_{n_{j}}\right| (2\mu _{0}-1)\pi }\right)^{n_{j}} 
              \label{s-esti1}
\end{eqnarray}

Similarly, 
\begin{equation}
\left| \frac{S_{n_{j}-1}(-n_{j}x_{n_{j}}(2\mu _{0}-1)\pi i)}
{(-ex_{n_{j}}(2\mu _{0}-1)\pi i)^{n_{j}}}\right| \leq 
K \left(\frac{e^{\varepsilon }}{\left| x_{n_{j}}\right| 
(2\mu _{0}-1)\pi } \right)^{n_{j}}  \label{s-esti2}
\end{equation}
By  (\ref{k-term})  and (\ref{g-inq})  the dominant term in 
$(ex_{n_{j}})^{-n_{j}}K_{n,\mu _{0}-1}(x_{n_{j}})$ 
corresponding to summation
indices $0\leq k\leq \mu _{0}-2$ is the term corresponding to $k=0$; that
is, the term 
\begin{equation*}
\frac{2g^{n_{j}}(x_{n_{j}})}{\pi i}
\left(1-\frac{1}{\sqrt{2\pi n_{j}}}
\frac{1}{1-x_{n_{j}}\pi i}g^{-n_{j}}(x_{n_{j}}) \left(1+O(n^{1-3\alpha }) \right) \right)
\end{equation*}
The other term with the same index $k=0$ is 
\begin{equation*}
\frac{2g^{n_{j}}(-x_{n_{j}})}{-\pi i}\left(1-\frac{1}{\sqrt{2\pi n_{j}}}
\frac{1}{1+x_{n_{j}}\pi i}g^{-n_{j}}(-x_{n_{j}}) \left(1+O(n^{1-3\alpha }) \right) \right)
\end{equation*}
which is still dominated by the above term. Note that 
\begin{equation*}
\left| g(x_{n_{j}})\right| =
\left| \frac{e^{\pi x_{n_{j}}i}}{e\pi x_{n_{j}}i}
\right| =\frac{e^{\pi \left| x_{n_{j}}\right| \sin \theta _{n_{j}}}}
{e\pi  \left| x_{n_{j}}\right| },
\end{equation*}
where $x_{n_{j}}=\left| x_{n_{j}}\right| e^{-i\theta _{n_{j}}}$.  Since 
$x_{n_{j}}\rightarrow x_{0}=r_{0}e^{-i\theta _{o}}$ with $0<\theta _{0}<\pi$, 
 this implies $\sin \theta _{n_{j}}\geq 0$ for all large $n_{j}$. Hence, 
$\left| g(x_{n_{j}})\right| \geq \frac{1}{e\pi \left| x_{n_{j}}\right| }$
implies 
\begin{equation}
\left| g(x_{n_{j}})\right| ^{n_{j}}\geq (\frac{1}{e\pi 
\left|  x_{n_{j}}\right| })^{n_{j}}.  \label{g-inq1}
\end{equation}
Note also that for $\mu _{0}\geq 2$ implies
$\frac{1}{e\pi }> \frac{e^{\varepsilon }}{(2_{\mu _{0}}-1)\pi }$.  
By comparing  (\ref{g-inq1}) with 
(\ref{s-esti1})  and (\ref{s-esti2})  we see the dominant term corresponding to $k=0$
still dominates the term corresponding to $k=\mu _{0}-1$.  Consequently, we
infer from equation  (\ref{m-kinzero1})  that the left hand side becomes
arbitrarily large as $n_{j}\rightarrow \infty$. Hence the left hand side
is non-zero, contradicting to the right side of the equation. This completes
the proof of part (b). 

We now proceed to prove part (c).  We show that if $x_{0}$ is a non-real
accumulation point in the region $\frac{1}{3\pi }<\left| x\right| <
\frac{1}{\pi }$, then either 
\begin{equation*}
\left| x_{0}\pi ie^{1-x_{0}\pi i}\right| =1
\quad
\textrm{or} 
\quad
\left| -x_{0}\pi ie^{1+x_{0}\pi i}\right| =1.
\end{equation*}
That $\frac{\pm 1}{\pi i}$ and real $x$ such that 
$\frac{-1}{\pi e}\leq \left| x\right| \leq \frac{1}{\pi }$ are points of the attractor is a
consequence of the density calculation. \ Since by assumption $x_{0}$ lies
in $\frac{1}{3\pi }<\left| x\right| <\frac{1}{\pi }$, we can certainly
choose $\varepsilon >0$ sufficiently small such that $x_{0}$ is in 
$\frac{1}{3\pi }+\varepsilon \leq \left| x\right| \leq \frac{1}{\pi }-\varepsilon$.
Again as in the previous cases we may assume 
\begin{equation*}
x_{0}=\left| x_{0}\right| e^{-i\theta _{0}},\text{where }0<\theta _{0}<\pi.
\end{equation*}
Also the same $\varepsilon $ works for an infinite sequence of zeros 
$x_{n_{j}}$ in 
\begin{equation*}
\frac{1}{3\pi }+\varepsilon \leq \left| x\right| \leq 
\frac{1}{\pi }-\varepsilon \text{ such that }x_{n_{j}}\rightarrow x_{0}\text{ as }
n_{j}\rightarrow \infty .
\end{equation*}
Now use Proposition  \ref{m-kprop}  with the choice $\mu =0$.
  So, 
\begin{equation*}
\frac{E_{n}(nx)}{n!}=M_{n,0}(x)+K_{n,0}(x),
\end{equation*}
where 
\begin{equation*}
M_{n,0}(x)=\frac{1}{\pi i}\oint_{\left| \xi \right| =1}
\left(\frac{e^{x\xi }}{\xi }\right)^{n}F_{0}(\xi )  \, d\xi ,
\end{equation*}
\begin{equation*}
K_{n,0}(x)=2 \left[\frac{S_{n-1}(nx\pi i)}{(\pi i)^{n+1}}+
\frac{S_{n-1}(-nx\pi i)}
{(-\pi i)^{n+1}} \right],
\end{equation*}
and 
\begin{equation}
F_{0}(\xi )=\frac{1}{\xi (e^{\xi }+1)}+
\frac{1}{\pi i(\xi -\pi i)}+
\frac{1}{(-\pi i)(\xi +\pi i)}.  \label{m-kinzero2}
\end{equation}
Now $E_{n_{j}}(n_{j}x_{n_{j}})=0$ implies 
$M_{n_{j,0}}(x_{n_{j}})+K_{n_{j,0}}(x_{n_{j}})=0$.  Since $F_{0}(\xi )$ is
analytic in the region $\varepsilon \leq \left| \xi \right| \leq 3\pi
-\varepsilon $, the asymptotics for $M_{n_{j,0}}(x_{n_{j}})$ is obtained as: 
\begin{equation}
M_{n_{j},0}(x_{n_{j}})=
\sqrt{\frac{2}{\pi }}  \left(   (x_{n_{j}}e)^{n_{j}}F_{0}
\left(     \frac{1}{x_{n_{j}}}\right)
\frac{1}{\sqrt{n_{j}}x_{n_{j}}} \right)   \,
\left(1+O(\frac{1}{n_{j}})    \right),
\label{m-esti}
\end{equation}
where the big $O$ term holds uniformly.  Note that since 
$x_{n_{j}}\rightarrow x_{0}$ we have $(\frac{1}{3\pi }+\varepsilon )\pi \leq
\left| x_{n_{j}}\pi i\right| \leq (\frac{1}{\pi }-\varepsilon )\pi $. So
$x_{n_{j}}\pi i$ lies in a compact set in the half plane $\Re  x<1$.
Hence, we can invoke Proposition  \ref{szego2} for the asymptotics for 
$K_{n_{j,0}}(x_{n_{j}})$. \ Combining this with (\ref{m-esti})  in
(\ref{m-kinzero2}) and expressing the result in terms of the function $g(x)$, 
we obtain
\begin{eqnarray*}
\lefteqn{
\sqrt{\frac{2}{\pi }}F_{0}(\frac{1}{x_{n_{j}}})
 \frac{1}{\sqrt{n_{j}}x_{n_{j}}}
\left(1+O(\frac{1}{n_{j}}) \right)   
}          \\
&& + \,
2\frac{g^{n_{j}}(x_{n_{j}})}{\pi i}-\sqrt{\frac{2}{\pi }}\frac{1}{\sqrt{n_{j}}
\left(1-x_{n_{j}} \right)\pi i} \left(1+O(n_{j}^{1-3\alpha }) \right) 
 \\
&& + \,
2\frac{g^{n_{j}}(-x_{n_{j}})}{(-\pi i)}-\sqrt{\frac{2}{\pi }}
\frac{1}{\sqrt{n_{j}}(1+x_{n_{j}})(-\pi i)} \left(1+O(n_{j}^{1-3\alpha })\right)=0.
\end{eqnarray*}
Note that there are terms in $\sqrt{\frac{2}{\pi }}F_{0}(1/x_{n_{j}})
\frac{1}{\sqrt{n_{j}}x_j}$ which are to be cancelled in the above equation and the
order term $O(\frac{1}{n_{j}})$ is absorbed in $O(n_{j}^{1-3\alpha })$.
This observation leads to a simplification: 
\begin{equation}
\sqrt{\frac{2}{\pi }}  \frac{1}  {e^{1/x_{n_{j}}         }+1    }
(1+O(n_{j}^{1-3\alpha }))  
    +
\sqrt{n_{j}}g^{n_{j}}(x_{n_{j}})
\left(\frac{2}{\pi i}+\frac{2}{(-\pi i)}
\left(\frac{g(-x_{n_{j}})}{g(x_{n_{j}})} \right)^{n_{j}} \right)=0  \label{g-eq}
\end{equation}

Also note that $\left| g(x_{0})\right| >\left| g(-x_{0})\right| $
 (strictly greater).  This implies 
\begin{equation*}
\left| \frac{g(-x_{n_{j}})}{g(x_{n_{j}})}\right| ^{n_{j}}
\rightarrow 0\text{ as }n_{j}\rightarrow \infty.
\end{equation*}
Therefore, if $\left| g(x_{0})\right| >1$,   then the left hand side of 
(\ref{g-eq})  becomes arbitrarily large in modulus for large $n_{j}$.
 This is a
contradiction. \ But if $\left| g(x_{0})\right| <1$,   then the left hand side
of (\ref{g-eq})   goes to $\sqrt{\frac{2}{\pi }}\frac{1}{e^{1/x_{0}}+1}$,  as 
$n_{j}\rightarrow \infty $, a non-zero number which is still a contradiction.
Hence we must have 
\begin{equation*}
\left| g(x_{0})\right| =1
\end{equation*}
that is, 
\begin{equation*}
\left| x_{0}\pi ie^{1-x_{0}\pi i}\right| =1.
\end{equation*}
This also shows that $z_{0}:=x_{0}\pi i$ is a point on the Szeg\"o curve: 
\begin{equation*}
\left| z_{0}\right| \leq 1,
\, 
\left| z_{0}e^{1-z_{0}}\right| =1,
\,
-\pi /2< \arg z_{0}<\pi /2.
\end{equation*}
that is, $\Re z_{0}>0$ or $z_{0}$ lies in the shaded segment of the
curve or $x_{0}$ is in the rotated Szeg\"o's curve.

To show there is no real accumulation point $x_{0}$ satisfying 
$\frac{1}{\pi } >x_{0}>\frac{1}{\pi e}$,  we can still use the asymptotics that leads to 
(\ref{g-eq}). In particular,  (\ref{g-eq})  still holds.  Note that in this case 
$x_{n_{j}}\pi i$ lies in the exterior region of the rotated Szeg\"{o} curve
in the $x$-plane. \ This implies $\left| g(x_{n_{j}})\right| <1$.  Hence in
the limit as $n_{j}\rightarrow \infty $ in  (\ref{g-eq})  we get 
\begin{equation*}
\sqrt{\frac{2}{\pi }}\frac{1}{e^{1/x_{0}}+1}=0,
\end{equation*}
a contradiction.  This finishes the proof of part (c).  
\end{proof}

This establishes the fact that the point set $K$ contains all points of the
zero attractor. 

\section{Determination of the density of zeros}

The next step is to carry out the density calculation.  The Euler
polynomials satisfy 
\begin{equation*}
E_{n}(1-x)=(-1)^{n+1}E_{n}(x)
\end{equation*}
We have

\begin{equation*}
E_{n}(nx)=(-1)^{n}E_{n}(1-nx)=(-1)^{n}E_{n}(-n(x-1/n)).
\end{equation*}
 From this it is easy to see that the zero attractor does has the
reflection symmetry with respect to the $y$-axis.  The Euler polynomials
are polynomials with real coefficients.  It obviously has the symmetry with
respect to reflection about the $x$-axis.  To describe the density we
choose the lower half of $K$ for discussion.

\begin{figure}[!ht]
\centerline{\psfig{figure=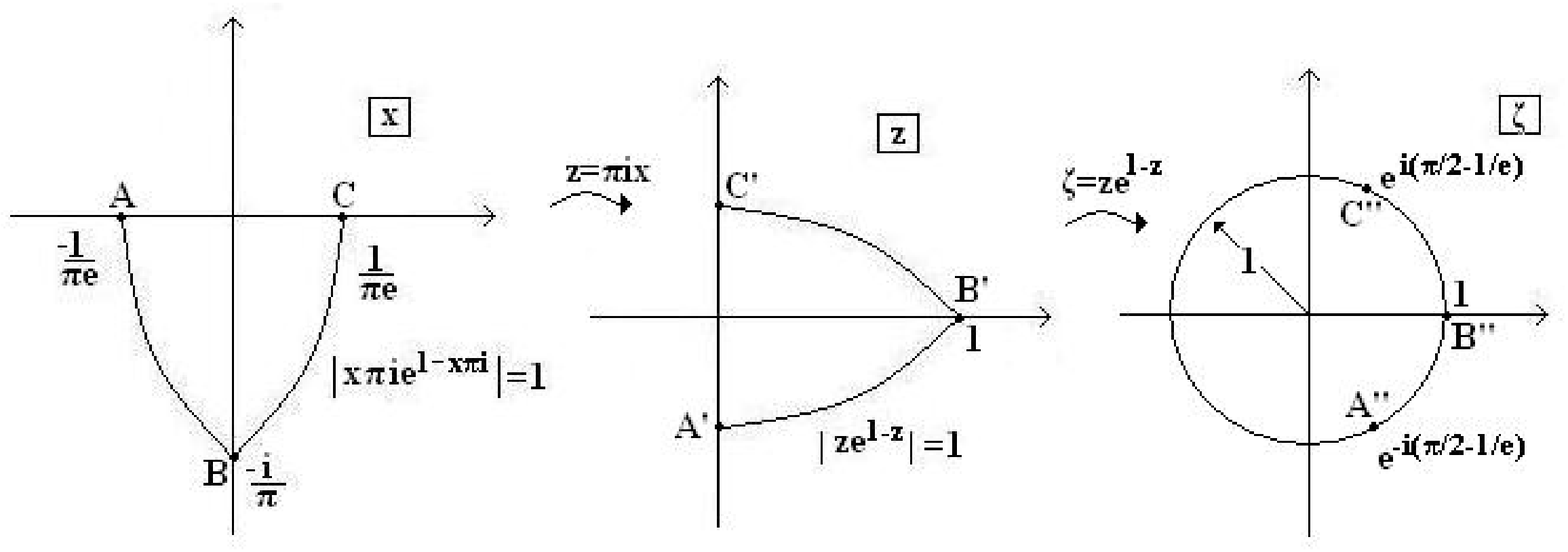,  width=4.0in}}
\caption{}
\label{fig:mappings_3}
\end{figure}

The image of the points $A,B,$ and $C$ in the $x$-plane under the mapping of 
$z=\pi ix$ are denoted by $A^{\prime }, B^{\prime }$,  and $C^{\prime }$
respectively. \ Further, the subsequent images $A^{\prime \prime }$, 
$B^{\prime \prime },$ and $C^{\prime \prime }$ are the images of $A^{\prime }$, 
$B^{\prime}$, and $C^{\prime }$ under $\zeta =ze^{1-z}$
(see Figure 3).

%
%

We show:

\begin{theorem}
\label{den-arc}
The image of the zeros of $E_{n}(nx)$ along the arc $BC$ in
the $x$-plane are uniformly distributed in the $\zeta $-plane along the
corresponding circular arc $B^{\prime \prime }C^{\prime \prime }$.  As a
consequence, the fraction of zeros residing in a neighborhood of the arc $BC$
in the $x$-plane is 
\begin{equation*}
\frac{\pi /2-1/e}{2\pi }=\frac{1}{4}-\frac{1}{2\pi e}.
\end{equation*}

\begin{theorem}
\label{den-seg}The real zeros of $E_{n}(nx)$ falling in the line segment 
$\overline{AC}$ are uniformly distributed in the segment $\overline{AC}$. 
As a consequence, the fraction of zeros residing in the segment $\overline{AC}$ is 
$\frac{2}{\pi e}$.
\end{theorem}
\end{theorem}

Before enumerating the zeros, we mention some well-known facts about the
Szeg\"{o} curve: The mapping $\zeta =ze^{1-z}$ is conformal  \cite{duren} in 
$\left| z\right| <1$.  The Szeg\"{o} curve is defined by:

\begin{equation*}
\left| ze^{1-z}\right| =1, \,  \left| z\right| \leq 1
\end{equation*}

The mapping $\zeta =ze^{1-z}$ is 
 in a neighborhood
of $z=1$. \ It is $2$-to-$1$ in a neighborhood of $z=1$ since 
$\frac{d\zeta }{dz}\left| _{z=1}=0\text{ and }
\frac{d^{2}\zeta }{dz^{2}}\right| _{z=1}\neq
0$.  Therefore, to enumerate the number of zero images inside a contour we
will face a difficulty of inverting the function in a neighborhood of $1$ in
the $\zeta $-plane. For this reason we will choose a contour  that does
not enclose $1$.  But how many zero images are left out in a neighborhood
of $1$?  We will show that it is of order $o(n)$.  To prove this we use
Jensen's inequality.  However, the use of Jensen's inequality requires a
knowledge of the function $E_{n}(nx)$ in a neighborhood of $-i/\pi$.
Fortunately, this knowledge is sufficiently provided by 
Proposition  \ref{szego1} and Proposition  \ref{szego2}.

\begin{theorem}
\label{jens}(Jensen's Inequality) \ Let $h(z)$ be analytic in the disc 
$\left| z-a\right| \leq R$ and let $0<r<R$ and $m$ be the number of zeros of 
$h(z)$ in the disc $\left| z-a\right| \leq r$. \ Then we have the following
inequality: 
\begin{equation*}
\left(  \frac{R}{r}   \right)^{m}
\leq 
\frac{   
\max_{  | z-a | =R}
\left|   h(z)  \right|     }
{        \left|   h(a)    \right|     }.
\end{equation*}
\end{theorem}

A proof can be found in most books of complex analysis. Let 
$\gamma_{\varepsilon }^{(n)}$ denote the number of zeros of $E_{n}(nx)$ which lie
in the disc $\left| x-\frac{1}{\pi i}\right| \leq \varepsilon $ in the $x$-plane.

\begin{proposition}
\label{zero-esti}
For all sufficiently small $\varepsilon >0$,  there exists 
$n_{0}(\varepsilon )$ such that for all $n\geq n_{0}$ we have 
\begin{equation*}
\gamma _{\varepsilon }^{(n)}\leq \frac{1}{\ln 2}
\left(n\ln (1+O(\varepsilon) \right)+\ln \frac{K}{\varepsilon }
\end{equation*}
\noindent
where the big $O$ constant and $K$ are absolute.
\end{proposition}

\begin{proof}
Since $\frac{\sqrt{n}}{(ex)^{n}}$ does not vanish in 
$\left| x-\frac{1}{\pi i}\right| \leq \varepsilon$,  the function 
\begin{equation}
h_{n}\left( x\right) :=\frac{E_{n}(nx)\sqrt{n}}{n!(ex)^{n}}  \label{h-def}
\end{equation}
is well-defined and has the same number of zeros as $E_{n}(nx)$ in 
$\left| x-\frac{1}{\pi i}\right| \leq \varepsilon$.
 We apply Theorem \ref{jens} to 
$h_{n}(x)$ on the disc 
$\left| x-\frac{1}{\pi i}\right| \leq 2\varepsilon$.
Thus 
\begin{equation*}
2^{\gamma _{\varepsilon }^{(n)}}=
\left(  \frac{2\varepsilon }{\varepsilon }  \right)^{\gamma _{\varepsilon }^{(n)}}
\leq \frac{\max_{\left| x-\frac{1}{\pi i}
\right| =2\varepsilon }\left| h_{n}(x)\right| }
{\left| h_{n}(\frac{1}{\pi i})\right| }.
\end{equation*}
We need asymptotic estimates for $\left| h_{n}(\frac{1}{\pi i})\right| $ and 
$\max_{\left| x-\frac{1}{\pi i}\right| =2\varepsilon }\left| h_{n}(x)\right|$.
 First, for $\left| h_{n}(\frac{1}{\pi i})\right| $ we apply Proposition
\ref{m-kprop}  with $\mu =0$.   Thus 
\begin{equation*}
\frac{E_{n}(nx)}{n!}=M_{n,0}(x)+K_{n,0}(x),
\end{equation*}
where 
\[
M_{n,0}(x)=\frac{2}{2\pi i}
\oint_{\left| \xi \right| =1}
\left(\frac{e^{x\xi }} {\xi } \right)^{n}
F_{0}(\xi ) \, d\xi ,
\]
and 
\[
K_{n,0}(x)=2
\left[\frac{S_{n-1}(nx\pi i)}{(\pi i)^{n+1}}+\frac{S_{n-1}(-nx\pi i)}
{(-\pi i)^{n+1}} \right].
\]
The asymptotics of $M_{n,0}(\frac{1}{\pi i})$ is known: 
\begin{equation}
M_{n,0}
\left(  \frac{1}{\pi i} \right)
=\sqrt{\frac{2}{\pi }  }
\left(    
          \left(  \frac{e}{\pi i}  \right)^{n}
F_{0}(   \pi i)  \frac{\pi i}{\sqrt{n}} 
\right)  \,
\left(1+O(\frac{1}{n})     \right).  \label{m-esti1}
\end{equation}

\noindent
Recall 
\begin{equation*}
K_{n,0}
\left(\frac{1}{\pi i} \right)=\frac{2}{(\pi i)^{n+1}}
\left[S_{n-1}(n)+\frac{S_{n-1}(-n)}{(-1)^{n+1}} \right].
\end{equation*}
We will show that $K_{n,0}(\frac{1}{\pi i})$ is the dominant term for 
$E_{n}(n\frac{1}{\pi i})$.  Use Proposition \ref{szego3} to get 
(with $t=n/(n-1)
$ and $n\rightarrow n-1)$:
\begin{equation}
S_{n-1}(n)=e^{n}
\left(\sqrt{\frac{2}{\pi }}\frac{\xi \cdot \frac{n}{n-1}}
{\frac{n}{n-1}-1} {\rm Erfc}(\sqrt{n-1}\xi )
\right) \,
\left(1+O(\frac{1}{\sqrt{n}}) \right),  \label{s-esti3}
\end{equation}
where $\xi =\left| \frac{n}{n-1}-1-\ln \frac{n}{n-1}\right| ^{1/2}=
\frac{1}{\sqrt{2}}\frac{1}{n-1} \left(1+O(\frac{1}{n}) \right)$.   So 
\begin{equation*}
\frac{\xi \frac{n}{n-1}}{\frac{n}{n-1}-1}=\frac{1}{\sqrt{2}}
\left(1+O(\frac{1}{n})\right).
\end{equation*}
Recall that 
\begin{eqnarray*}
{\rm Erfc}(x) &=&
\int_{x}^{\infty }e^{-s^{2}}ds \\
&=&
\int_{0}^{\infty }e^{-s^{2}}ds-\int_{0}^{x}e^{-s^{2}}ds \\
&=&
\frac{\sqrt{\pi }}{2}+O(x)
\end{eqnarray*}
as $x\rightarrow 0$. 
Hence 
\begin{eqnarray*}
{\rm Erfc}(\sqrt{n-1}\xi )
&=&
{\rm Erfc}\left(\frac{1}{\sqrt{2}}\frac{1}{\sqrt{n-1}}
(1+O(\frac{1}{n})) \right)\\
&=&
\frac{\sqrt{\pi }}{2}+O(\frac{1}{\sqrt{n}}).
\end{eqnarray*}
Inserting these estimates back into  (\ref{s-esti3})  we have 
\[
S_{n-1}(n)=\frac{e^{n}}{2}
\left(1+O(\frac{1}{\sqrt{n}}) \right).
\]
Similarly using Proposition \ref{szego3} we have 
\[
S_{n-1}(-n)=O(e^{-n}).
\]
These estimates give 
\[
K_{n,0}(\frac{1}{\pi i})=
(\frac{1}{\pi i})(\frac{e}{\pi i})^{n}
\left(1+O(\frac{1}{\sqrt{n}}) \right).
\]
Comparing (\ref{m-esti1})  with the above we see that the order of 
$M_{n,0}(\frac{1}{\pi i})$ is small than that of $K_{n,0}(\frac{1}{\pi i})$ by a
factor of $\sqrt{n}$.   Hence 
\begin{equation*}
\frac{E_{n}(n\frac{1}{\pi i})}{n!}=(\frac{1}{\pi i})
(\frac{e}{\pi i})^{n} \left(1+O(\frac{1}{\sqrt{n}})\right),
\end{equation*}
and 
\[
h_{n}(\frac{1}{\pi i}):=\frac{E_{n}(n\frac{1}{\pi i})\sqrt{n}}{n!(\frac{e}
{\pi i})^{n}}=\frac{\sqrt{n}}{\pi i}  \left(1+O(\frac{1}{\sqrt{n}}) \right).
\]
We now estimate $\max_{\left| x-\frac{1}{\pi i}\right| =2\varepsilon }\left|
h_{n}(x)\right|$.  To this end we use Proposition \ref{m-kprop} with $\mu =0$:
\[
\frac{E_{n}(nx)}{n!}=M_{n,0}(x)+K_{n,0}(x),
\]
where 
\begin{eqnarray*}
M_{n,0}(x)
&=&
\frac{1}{\pi i}  \oint_{\left| \xi \right| =1}
\left(\frac{e^{x\xi }}{\xi } \right)^{n}F_{0}(\xi )  \, d\xi ,
\\
K_{n,0}(x) 
&=&
2 \left[\frac{S_{n-1}(nx\pi i)}{(\pi i)^{n+1}}+
\frac{S_{n-1}(-n)}{(-\pi  i)^{n+1}} \right].
\end{eqnarray*}
The asymptotics of $M_{n,0}(x)$ is: 
\[
M_{n,0}(x)=
\sqrt{\frac{2}{\pi }}
\left(  (xe)^{n}F_{0}(\frac{1}{x})
\frac{1}{\sqrt{n}x}  \right)  \,  \left(1+O(\frac{1}{n})\right).
\]
Care must be exercise to handle the asymptotics of $S_{n-1}(nx\pi i)$.
Recall that $x$ is on the circumference 
$\left| x-\frac{1}{\pi i}\right| =2\varepsilon $.

\begin{figure}[ht]
\centerline{\psfig{figure=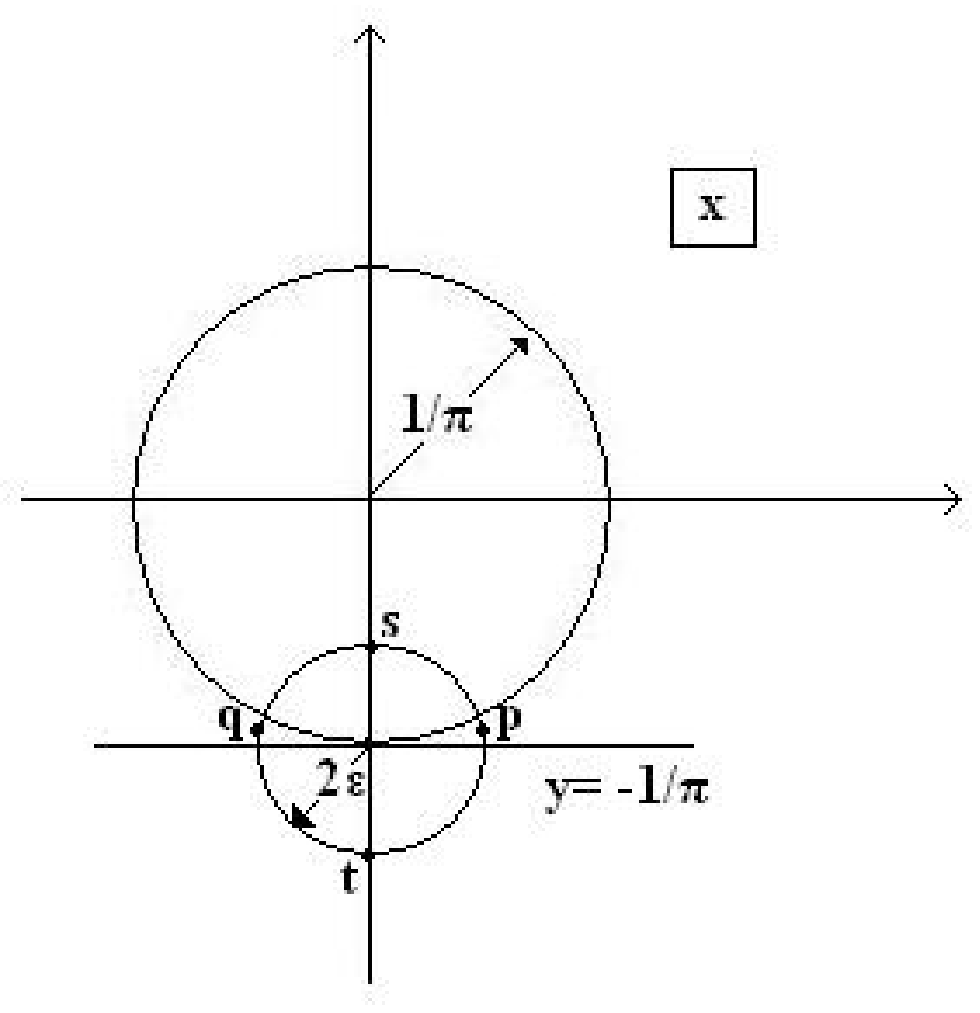, width=3.5in}}
\caption{}
\label{fig:circle_4}
\end{figure}

Choose two points $p$ and $q$ (see Figure 4)
such that $p$ is
the mid-point of the circular arc between the circle 
$\left| x\right| =  \frac{1}{\pi }$ and the horizontal tangent line to the circle 
$\left| x\right| = \frac{1}{\pi }$at the point $x=\frac{1}{\pi i}$.
 The point $q$ is similarly
selected. Note that as $\varepsilon \rightarrow 0^{+}$, 
 the distance between  $p$ and the circle $\left| x\right| =\frac{1}{\pi }$ is of order 
$O(\varepsilon ^{2})$. This is so because the horizontal line 
$y=\frac{-1}{\pi }$ is tangent to the circle $\left| x\right| =\frac{1}{\pi }$.  The image
of the arc $qtp$ under the map $z=x\pi i$ is a circular arc in the 
$z$-plane where Proposition 1 is applicable. Thus for all $x$ on the arc $qtp$
we have 
\begin{equation}
S_{n-1}(nx\pi i)=\frac{e^{nx\pi i}g^{-n}(x)}{\sqrt{2\pi n}(x\pi i-1)}
\left(1+O_{\varepsilon }(n^{(1-3\alpha )}) \right),  \label{s-esti4}
\end{equation}
where $O_{\varepsilon }$ stands for the $\varepsilon $-dependence for the
big $O$ constant. \ Taking  (\ref{s-esti4})  into consideration we get 
\begin{eqnarray*}
\lefteqn
{
h_{n}(x)=\sqrt{\frac{2}{\pi }}F_{0}(\frac{1}{x})\frac{1}{x}
(1+O(\frac{1}{n}))
} \\
&& \quad
+ \,
\frac{2g^{n}(x)}{\pi i}\cdot \frac{g^{-n}(x)}{\sqrt{2\pi }(x\pi i-1)}
(1+O_{\varepsilon }(n^{1-3\alpha }))  \\
&& \quad
+ \,
\frac{2\sqrt{n}g^{n}(-x)}{(-\pi i)}
\left[1-\frac{g^{-n}(-x)}{\sqrt{2\pi n}
(1+x\pi  i)}(1+O_{\varepsilon }(n^{1-3\alpha }))  \right],
\end{eqnarray*}
where the third term is obtained from an analogous application of Proposition 
\ref{szego2} to $S_{n-1}(-nx\pi i)$. We now estimate $h_{n}(x)$ as follows.
The first term 
\[
\left| \sqrt{\frac{2}{\pi }}F_{0}(\frac{1}{x})\frac{1}{x}
\left(1+O(\frac{1}{n}) \right)\right| 
\]
 is obviously $\leq K$, an absolute constant. The
second term 
\begin{eqnarray*}
\left| \frac{2g^{n}(x)}{\pi i}\frac{g^{-n}(x)}{\sqrt{2\pi }(x\pi i-1)}
(1+O_{\varepsilon }(n^{1-3\alpha }))\right| 
&\leq &
\frac{K}{\left| x\pi i-1\right| }(1+O_{\varepsilon }(n^{1-3\alpha })) \\
&\leq& 
\frac{K}{\varepsilon }(1+O_{\varepsilon }(n^{1-3\alpha }))
\end{eqnarray*}
The third term 
\begin{equation*}
\left| \frac{2\sqrt{n}g^{n}(-x)}{(-\pi i)}
\left[1-\frac{g^{-n}(-x)}{\sqrt{2\pi n}
(1+x\pi i)}(1+O_{\varepsilon }(n^{1-3\alpha })) \right]\right| \leq K,
\end{equation*}
because $-x\pi i$ lies outside of the Szeg\"{o} curve and 
$\left| g(-x)\right| <1$. \ Hence for $x$ on the arc $ptq$ on 
$\left| x-\frac{1}{\pi  i}\right| =2\varepsilon $ we have 
\begin{equation*}
\left| h_{n}(x)\right| \leq \frac{K}{\varepsilon }
\left(1+O_{\varepsilon }(n^{1-3\alpha }) \right),
\end{equation*}
where $K$ is an absolute constant. When $x$ in the arc $psq$ on the
circumference$\left| x-\frac{1}{\pi i}\right| =2\varepsilon $, we invoke
Proposition \ref{szego2}. 
 Thus 
\begin{eqnarray}
\lefteqn{
h_{n}(x)=\sqrt{\frac{2}{\pi }} \left(F_{0}(\frac{1}{x})\frac{1}{x}  \right)
\left(1+O(\frac{1}{n}) \right)+
\frac{2\sqrt{n}g^{n}(x)}{\pi i}
} \nonumber \\
&& \qquad \times \,
\left[1-\frac{g^{-n}(x)}{\sqrt{2\pi n}(1-x\pi i)}
(1+O_{\varepsilon }(n^{1-3\alpha })) \right]
\nonumber \\
&&
\qquad + \,
\frac{2\sqrt{n}g^{n}(-x)}{(-\pi i)}
\left[1-\frac{g^{-n}(-x)}{\sqrt{2\pi n}
(1+x\pi i)}(1+O_{\varepsilon }(n^{1-3\alpha })) \right].  \label{h-esti}
\end{eqnarray}
The magnitude of $g(x)$ for $x$ on $psq$ is estimated as follows.  Recall 
$ x\pi i=1+2\varepsilon \pi ie^{i\theta }$.  This implies 
\begin{eqnarray*}
\lefteqn
{
\left| x\pi ie^{1-x\pi i}\right| =\left| x\pi \right| \left|
e^{-2\varepsilon \pi e^{i(\theta +\pi /2)}}\right| 
}\\
&&
\qquad
=\left| x\pi \right| e^{-2\varepsilon \pi \cos (\theta +\pi /2)}
\geq \left| 1-2\varepsilon \pi \right| e^{-2\varepsilon \pi }
\end{eqnarray*}
This yields 
\begin{equation*}
\left| g(x)\right| =\left| \frac{1}{x\pi ie^{1-x\pi i}}\right| \leq 
\frac{
e^{2\varepsilon \pi }}{1-2\varepsilon \pi }.
\end{equation*}
Therefore, we infer from  (\ref{h-esti})  that 
\begin{equation*}
\left| h_{n}(x)\right| \leq 
K\left[
\left(\frac{e^{2\varepsilon \pi }}{1-2\varepsilon \pi } \right)^{n}\sqrt{n}
+O_{\varepsilon }(n^{1-3\alpha }) \right]
\end{equation*}
Note that in deriving the above equation the fact that $\left| g(-x)\right|
<1$ was still used. Combine these estimates to get 
\begin{equation*}
\max_{\left| x-\frac{1}{\pi i}\right| =2\varepsilon }\left| h_{n}(x)\right|
\leq K\left[(\frac{e^{2\varepsilon \pi }}{1-2\varepsilon \pi })^{n}\sqrt{n}
+O_{\varepsilon }(n^{1-3\alpha }) \right].
\end{equation*}
Recall that 
\begin{eqnarray*}
2^{\gamma _{\varepsilon }^{(n)}}
&\leq& \frac{\max_{\left| x-\frac{1}{\pi i}\right| 
=
2\varepsilon }\left| h_{n}(x)\right| }{\left| h_{n}(\frac{1}{\pi i})\right| }
\\
&\leq&
 \frac{K[(\frac{e^{2\varepsilon \pi }}{1-2\varepsilon \pi })^{n}\sqrt{n}
+O_{\varepsilon }(n^{1-3\alpha })]}{\frac{\sqrt{n}}{\pi }
\left(1+O(\frac{1}{\sqrt{n}})  \right)}
\end{eqnarray*}
This implies that 
\begin{equation*}
\gamma _{\varepsilon }^{(n)}
\leq \ln  
\left[K \left(\frac{e^{2\varepsilon \pi }}
{1-2\varepsilon \pi } \right)^{n}+
O_{\varepsilon }(n^{1/2-3\alpha }))
\right]/\ln 2.
\end{equation*}
This finishes the proof of  Proposition \ref{zero-esti}.  
\end{proof}

We now come to determine the density of zeros.  First of all, we refer to
the mapping relation in Figure 3.  In general, we let $N_{n}(\alpha ,\beta )$
be the number of image points of zeros of $E_{n}(nx)$ in the $\zeta$-plane
that fall in the angular sector $\alpha \leq \arg \zeta \leq \beta$. Now
let an arbitrary $\theta$ be given in the interval

\begin{equation*}
0<\theta <\frac{\pi }{2}-\frac{1}{e}.
\end{equation*}

Our goal is to establish the following proposition:

\begin{proposition}
\label{den-arc-pr} 
\quad $\displaystyle
\lim_{n\rightarrow \infty }\frac{1}{n}
N(0,\theta )=\frac{\theta }{2\pi }.
$
\end{proposition}

\begin{proof}
Recall that $N_{n}(0,\theta )$ is the number of image points in the $\zeta $
-plane that fall in the angular sector $0\leq \arg \zeta \leq \theta$. We
alleviate the problem that the straight edges of the above sector may
contain some image points by perturbation.  The reason for requiring that
no image points of zeros fall on the straight edges is to guarantee an
application of the argument principle.  Since the Euler polynomials are
polynomials of rational coefficients, the roots are algebraic numbers in the 
$x$-plane.  This implies that the totality of image points in the 
$\zeta$-plane is countable.  Hence we can choose an arbitrary small number 
$\varepsilon >0$ so that the straight edges of $3\varepsilon =\arg \zeta $
avoids all images points of zeros.  The following arguments assume a fixed 
$\varepsilon >0$ and $n$ will be $\rightarrow \infty$ eventually.  It is
now obvious that we have the following inequality: 
\begin{equation}
N_{n}(0,\theta )-N_{n}(3\varepsilon ,\theta +
\varepsilon _{n}^{\prime})=
O \left(\gamma _{\varepsilon }^{(n)} \right),  \label{n-esti}
\end{equation}
where a null sequence $\varepsilon _{n}^{\prime }>0$ is chosen so that no
image points are on the straight edge 
$\arg \zeta =\theta +\varepsilon_{n}^{\prime }$.  Define 
\begin{equation*}
f_{n}(\zeta ):=\frac{E_{n}(nx(\zeta ))}{n!(ex(\zeta ))^{n}}\sqrt{n}.
\end{equation*}
Note that:  (1)  the image points are zeros of $f_{n}(\zeta )$
and
(2)
$f_{n}(\zeta)=h_{n}(x(\zeta ))$, where $h_{n}(x)$ was defined in  (\ref{h-def}). 
Here the
function $x(\zeta )$ is the inverse map of the map from the $x$-plane to the 
$\zeta $-plane. \ $x(\zeta )$ is 1-1 except in a small neighborhood of 
$\zeta =1$.  

\begin{figure}[ht]
\centerline{\psfig{figure=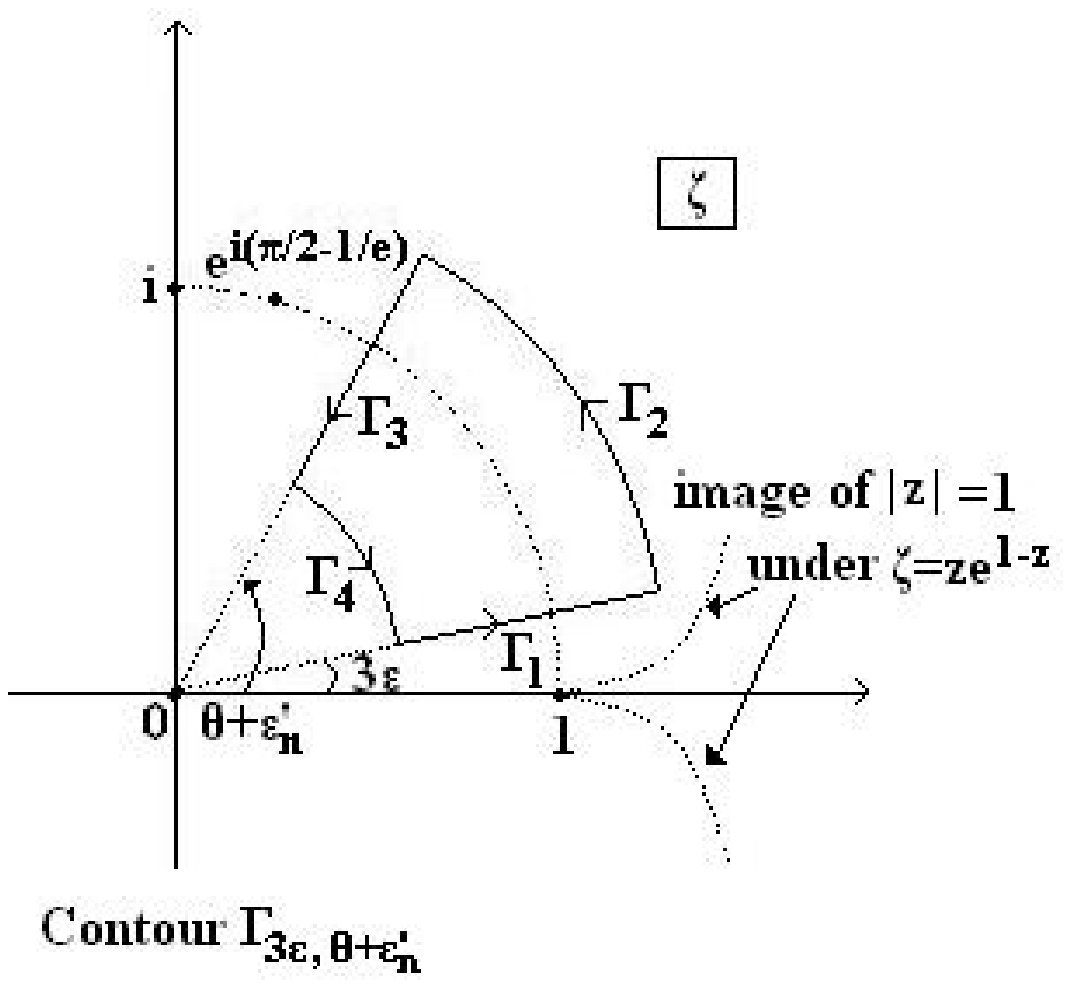, height=3.0in, width=4.0in}}
\caption{}
\label{fig:contour_5}
\end{figure}

%
%

We introduce the contour 
$\Gamma _{3\varepsilon ,\theta +\varepsilon_{n}^{\prime }}$ defined as
(see Figure 5):
$\Gamma _{1}\cup \Gamma _{2}\cup \Gamma_{3}\cup \Gamma _{4}$:
\[
\left\{
\begin{array}{lll}
\Gamma _{1}&:&  \, re^{i3\varepsilon },
 1-\varepsilon /2\leq r\leq 1+\varepsilon /2,
\\
\Gamma _{2} &:&  \, (1+\varepsilon /2)e^{i\phi },
   3\varepsilon \leq \phi \leq \theta
+\varepsilon _{n}^{\prime },
\\
\Gamma _{3} &:&  \, re^{i(\theta +\varepsilon _{n}^{\prime })},
1-\varepsilon /2\leq  r\leq 1+\varepsilon /2,
\\
\Gamma _{4} &:&  \, (1-\varepsilon /2)e^{i\phi },
3\varepsilon \leq \phi \leq \theta
+\varepsilon _{n}^{\prime }.
\end{array}
\right.
\]
We apply the argument principle to enumerate the number of image points. 
Thus, 
\begin{eqnarray}
N_{n}(3\varepsilon ,\theta +\varepsilon _{n}^{\prime }) 
&=&  \frac{1}{2\pi i}
\oint_{\Gamma _{3\varepsilon ,\theta +\varepsilon _{n}^{\prime }}}
\frac{  \frac{d}{d\zeta }f_{n}(\zeta )}{f_{n}(\zeta )}  \,  d\zeta 
\nonumber \\
&=&\frac{1}{2\pi }
\Im(\oint_{\Gamma _{3\varepsilon ,\theta +
\varepsilon_{n}^{\prime }}}
\frac{\frac{d}{d\zeta }f_{n}(\zeta )}{f_{n}(\zeta )}  \, d\zeta ) 
\nonumber  \\
&=&
\frac{1}{2\pi } \left(\Im(I_{1})+\Im(I_{2})+\Im(I_{3})+
\Im (I_{4})  \right),  \label{4integ}
\end{eqnarray}
where $I_{i}=\int_{\Gamma _{i}}, \, 1\leq i  \leq 4$.
We shall handle $I_{2}$
first. If we pull the contour $\Gamma _{2}$ back in the $z$-plane, it falls
in a compact set in the half plane $\Re(z)<1$ so that we can use
Propsition  \ref{szego2} for the asymptotics of $S_{n-1}(nz)$.  Thus with an
application of Proposition \ref{m-kprop} (with $\mu =0$) we have 
\begin{eqnarray*}
\lefteqn{
\frac{E_{n}(nx) \sqrt{n}}{n!(ex)^{n}}
=
\sqrt{\frac{2}{\pi }}F_{0}(\frac{1}{x})
\frac{1}{x}  \left(1+O(\frac{1}{n}) \right)
}
\\
&& \qquad \quad
+ \,
\frac{2\sqrt{n}  g^{n}(x)}  {\pi i}      
\left[1-\frac{g^{-n}(x)}{\sqrt{2\pi n}(1-x\pi i)}
\left(1+O_{\varepsilon }(n^{1-3\alpha }) \right) \right]
\\
&& \qquad \quad
+ \,
\frac{2\sqrt{n}g^{n}(-x)}      {(-\pi i) }            
\left[
1-\frac{  g^{-n}(-x)      }  {      \sqrt{2\pi n}         (1+x\pi i)  }                  
  (1+O_{\varepsilon }(n^{1-3\alpha }))          
\right]
\end{eqnarray*}
Recall 
\begin{equation*}
F_{0}(\xi )=\frac{1}{\xi (e^{\xi }+1)}+\frac{1}{\pi i(\xi -\pi i)}+
\frac{1}{ (-\pi i)(\xi +\pi i)},
\end{equation*}
so 
\begin{equation*}
F_{0}(\frac{1}{x})=\frac{x}{e^{1/x}+1}+\frac{1}{\pi i(1-x\pi i)}+
\frac{1}{(-\pi i)(1+x\pi i)}.
\end{equation*}
 Canceling some common terms to simplify the expression we get 
\begin{eqnarray}
\lefteqn{
\frac{E_{n}(nx)\sqrt{n}}   {n!(ex)^{n}}
=  \sqrt{\frac{2}{\pi }}      \frac{1}{e^{1/x}+1} \, \times
} \nonumber  \\
&&
\,
\left\{  1+O(\frac{1}{n} ) 
+
\frac{g^{n}(x)\sqrt{2\pi n}(e^{1/x}+1)}{\pi i} 
+
\frac{g^{n}(-x)\sqrt{2\pi n}(e^{1/x}+1)}{(-\pi i)}+O(n^{1-3\alpha })  \right\}.
\label{h-def1}
\end{eqnarray}

Furthermore since the pulled back $x$'s lie outside the zero attractor so 
$\left| g(-x)\right| <\left| g(x)\right| <1$.  In this way we see that 
\begin{equation*}
\lim_{n\rightarrow \infty }\frac{E_{n}(nx)\sqrt{n}}{n!(ex)^{n}}
\frac{e^{1/x}+1}{\sqrt{2/\pi }}\rightarrow 1, \, \text{ uniformly}
\end{equation*}
or 
\begin{equation*}
h_{n}(x)\rightarrow \frac{\sqrt{2/\pi }}{e^{1/x}+1},\, \text{ uniformly}
\end{equation*}
or equivalently 
\begin{equation*}
h_{n}(x(\zeta ))\rightarrow \frac{\sqrt{2/\pi }}{e^{1/x(\zeta )}+1}, \,
\text{  uniformly for }\zeta \in \Gamma _{2}.
\end{equation*}
By a theorem of uniform convergence of analytic functions we get 
\begin{equation*}
\frac{\frac{d}{d\zeta }h_{n}(x(\zeta ))}{h_{n}(x(\zeta ))}\rightarrow 
\frac{e^{1/x(\zeta )}\frac{dx/d\zeta }{x^{2}}}{e^{1/x(\zeta )}+1}
\end{equation*}
uniformly for $\zeta \in \Gamma _{2}$. 
Recall 
$f_{n}(\zeta )=h_{n}(x(\zeta ))$. 
Hence by integrating the above along $\Gamma _{2}$ we get 
\begin{eqnarray*}
\frac{1}{2\pi }\Im(I_{2})
&=&
\frac{1}{2\pi }\Im
\left(
\int_{\Gamma _{2}}
\frac{\frac{d}{d\zeta }f_{n}(\zeta )}{f_{n}(\zeta )}  \, d\zeta \right)   \\
&\rightarrow &
\frac{1}{2\pi }\Im
\left(
\int_{\Gamma _{2}}\frac{e^{1/x(\zeta )}
\frac{dx/d\zeta }{x^{2}}}{e^{1/x(\zeta )}+1}  
\, d\zeta \right).
\end{eqnarray*}
This implies 
\begin{equation*}
\frac{1}{n}\frac{1}{2\pi }\Im(I_{2})=O(\frac{1}{n}).
\end{equation*}
Let us handle $I_{4}$ now. The asymptotics in  (\ref{h-def1})  still hold good
in this case. But since $x$ now is inside the zero attractor so 
\begin{equation*}
\left| g(x)\right| >1\text{ and still }\left| g(x)\right| >
\left| g(-x)\right| .
\end{equation*}
Rewrite  (\ref{h-def1})  in the form 
\begin{eqnarray}
\lefteqn
{
\frac{f_{n}(\zeta )}{\sqrt{2\pi n}}\zeta ^{n}=\frac{1}{\sqrt{2\pi n}}
\frac{\zeta ^{n}\sqrt{\frac{2}{\pi }}}{e^{1/x}+1}(1+O(\frac{1}{n})+
O_{\varepsilon}(n^{1-3\alpha }))
} \nonumber \\
&& \qquad
+
\frac{\sqrt{2/\pi }}{\pi i}+\frac{\sqrt{2/\pi }}{(-\pi i)}(\zeta g(-x))^{n}.
\label{h-def2}
\end{eqnarray}
Note that $\zeta =ze^{1-z}$ and $\left| \zeta \right| =1-\varepsilon /2$
since $\zeta \in \Gamma _{4}$.  Also, 
$\left| \zeta g(-x)\right| =\left| g(-x)/g(x)\right| <1$. Hence we see that 
\[
\lim_{n\rightarrow \infty }
\frac{f_{n}(\zeta )}{\sqrt{2\pi n}}\zeta^{n}\rightarrow \frac{\sqrt{2/\pi }}{\pi i}
\]
uniformly and we get similarly 
\[
\frac{\frac{d}{d\zeta }f_{n}(\zeta )}{f_{n}(\zeta )}+\frac{n}{\zeta }
\rightarrow 0
\]
uniformly for $\zeta \in \Gamma _{4}$. 
Integration along $\Gamma _{4}$ gives 
\[
\Im(I_{4})=-\Im
\left(\int_{\Gamma _{4}}\frac{\frac{d}{d\zeta }
f_{n}(\zeta )}{f_{n}(\zeta )} \, d\zeta \right)
=-n\Im
\left(\int_{\Gamma _{4}}
\frac{ d\zeta }{\zeta } \right)+o(1).
\]
Note that the orientation of $\Gamma _{4}$ gives the correction of the 
``{$-$}''   sign.  Hence, 
\begin{equation*}
\Im(I_{4})=n
\Im \left(i(\theta +\varepsilon _{n}^{\prime}-3\varepsilon ) \right)+o(1).
\end{equation*}
This implies that 
\[
\frac{1}{n}\frac{1}{2\pi }\Im(I_{4})=\frac{1}{2\pi }
\left((\theta+\varepsilon _{n}^{\prime }-3\varepsilon ) \right)+o(\frac{1}{n}).
\]
Caution must be exercised when we handle $I_{1}$
( $I_{3}$ can be analogously
handled).  Although  (\ref{h-def2})  still holds good for $I_{1}$, but 
$\left| \zeta \right| $ varies from $1-\varepsilon /2$ to $1+\varepsilon /2$.   It
is not clear how  (\ref{h-def2})  becomes useful to us. We now employ a
number-theoretic argument and Theorem  \ref{jens} to obtain a useful estimate
for $I_{1}$. The argument we follow is inspired by that of K. Chandrasekharan
\cite{chandra}. Let $l$ be the number of points on $\Gamma _{1}$ so that 
$\Re(f_{n}(\zeta ))=0$.  We insert these points into $\Gamma _{1}$ and
decompose $I_{1}$ as the sum of integrals whose end points are two
\textit{consecutive}   points where $\Re(f_{n}(\zeta ))=0$ except possibly the
beginning integral and the final integral.  A typical integral to consider
is of the form $\Im(\int_{a}^{b}\frac{f_{n}^{\prime }(\zeta )}
{f_{n}(\zeta )} \,  d\zeta )$ , where $a$ and $b$ are two consecutive roots of 
$\Re(f_{n}(\zeta ))=0$ along $\Gamma _{1}$.  By a change of variable,
we have 
\begin{equation*}
\Im
      \left(
\int_{a}^{b}
\frac{  f_{n}^{  \prime }(  \zeta )  }  {  f_n (\zeta ) } 
\, d \zeta  
     \right)
=
\Im
\left(
\int_{C_n }   \frac{d \xi }  {\xi } 
   \right),
\end{equation*}
where the contour $C_{n}$ is the image of the segment $\overline{ab}$ under
the map $\zeta \rightarrow f_{n}(\zeta )$. We comment that $C_{n}$ 
intersects the imaginary axis at its two endpoints and at no other interior
point of $C_n$.
No matter what, we
always have by Cauchy's Theorem
\[
\int_{C_{n}}\frac{d\xi }{\xi }=\int_{K_{ab}}\frac{d\xi }{\xi },
\]
where $K_{ab}$ denotes the semi-circle with segment 
$\overline{f_{n}(a)f_{n}(b)}$ as the base so that $K_{ab}$ lies in the same half plane
as $C_{n}$.  Note that $\int_{K_{ab}}d\xi /\xi $ is either $0$,$ i \pi$,  or 
$-i\pi$.  Hence 
\[
\left| \Im\int_{a}^{b}\frac{f_{n}^{\prime }(\zeta )}{f_{n}(\zeta )} \,
d\zeta \right| =\left| \Im\int_{K_{ab}}\frac{d\xi }{\xi }\right| \leq
\pi .
\]
It follows that 
\begin{equation}
\left| \Im(I_{1})\right| \leq (l+1)\pi .  \label{l-ineq2}
\end{equation}
Here we have $(l+1)\pi $ as an upper estimate because of the possible
inclusion of the beginning and final integrals.  To estimate $l$,  we shall
use the Jensen's  inequality:  Define 
\[
\widetilde{h}_{n}(\zeta ):=
\frac{h_{n}(x(\zeta e^{i3\varepsilon}))+h_{n}(-x(\zeta e^{-i3\varepsilon }))}{2},
\]
where $h_{n}(x)$ was defined in  (\ref{h-def}). 
Let us note the following two
properties of $\widetilde{h}_{n}(\zeta )$:
\begin{enumerate}
\item $\widetilde{h}_{n}(\zeta )$ is
an analytic function of $\zeta$. 
\item When $\zeta $ is real and 
$1-\varepsilon/2\leq \zeta \leq 1+\varepsilon /2$,
 then $\zeta e^{i3\varepsilon } \in
\Gamma _{1}$ and $x(\zeta e^{ie\varepsilon })$ and further
$-x(\zeta e^{-ie\varepsilon })$ are complex conjugates in the $x$-plane. 
\end{enumerate}
\noindent
Hence 
\[
\widetilde{h}_{n}(\zeta )=
\Re
\left[h_{n}(x(\zeta e^{i3\varepsilon })) \right].
\]
The reason we choose $3\varepsilon$ instead of $\varepsilon $ is for
convenience as the sequel will show.  Now we regard 
$\widetilde{h}_{n}(\zeta )$
as an analytic function defined on the disc 
$\left| \zeta -(1+\varepsilon /2)\right| \leq 2\varepsilon $ in the $\zeta $-plane. 
(This is so because the circle 
$\{\zeta e^{i3\varepsilon }:\left| \zeta -(1+\varepsilon /2)\right| =2\varepsilon \}$ 
does not include $\zeta =1$ in its interior,
for otherwise, $x(\zeta e^{i3\varepsilon })$ would not be well-defined in
the 
$\{\zeta e^{i3\varepsilon }:\left| \zeta -(1+\varepsilon /2)\right| \leq
2\varepsilon \}$.) 
Note that the disc 
$\{\zeta :\left| \zeta -(1+\varepsilon /2)\right| \leq \varepsilon \}$ 
contains the real interval 
$[1-\varepsilon/2,1+\varepsilon /2]$. Each root of $\Re(f_{n}(\zeta ))$ along 
$\Gamma_{1}$ is a real root of $\widetilde{h}_{n}(\zeta )$ in 
$[1-\varepsilon/2,1+\varepsilon /2]$.
Let $\widetilde{l}$ be the number of roots of 
$\widetilde{h}_{n}(\zeta )$ in 
$\{\zeta :\left| \zeta -(1+\varepsilon /2)\right| \leq \varepsilon \}$
 (possible complex roots are counted in 
$\widetilde{l})$,  then obviously 
\begin{equation}
l\leq \widetilde{l}.  \label{l-ineq}
\end{equation}
Apply Theorem \ref{jens} in the disc 
$\{\zeta :\left| \zeta -(1+\varepsilon/2)\right| \leq 2\varepsilon \}$ 
to get 
\begin{equation}
\left( \frac{2\varepsilon }{\varepsilon }\right) ^{\widetilde{l}}
=2^{\widetilde{l}}\leq 
\frac{\max_{\left| \zeta -(1+\varepsilon /2)\right|
=2\varepsilon }\left| \widetilde{h}_{n}(\zeta )\right| }
{\left| \widetilde{h}_{n}   (1+\varepsilon /2)
\right| }.  \label{l-ineq1}
\end{equation}
Recall the definition of $\widetilde{h}_{n}(1+\varepsilon /2)$: 
\[
\widetilde{h}_{n}(1+\varepsilon /2)=\frac{1}{2}
\left[h_{n}(x((1+\varepsilon /2)e^{i3\varepsilon }))
+h_{n}(-x((1+\varepsilon /2)e^{-i3\varepsilon }))    \right]
\]
The point $\zeta =(1+\varepsilon /2)e^{3i\varepsilon }$ lies in a region in
the $\zeta $-plane so that the pull back 
$z((1+\varepsilon/2)e^{3i\varepsilon })$ lies in the half plane $\Re(z)<1$.
Therefore, we can use  (\ref{h-esti})  to determine the asymptotics of 
$h_{n}(x((1+\varepsilon /2)e^{3i\varepsilon }))$.  Similarly for 
$h_{n}(-x((1+\varepsilon /2)e^{-3i\varepsilon }))$. In this case we have 
\[
\left| g(-x((1+\varepsilon /)e^{-3i\varepsilon }))\right| 
<\left| g(x((1+\varepsilon /2)e^{3i\varepsilon }))\right| <1,
\]
so 
\begin{equation}
\lim_{n\rightarrow \infty }\widetilde{h}_{n}(1+\varepsilon /2)
=
\sqrt{\frac{2}{\pi }}\Re
\left[\frac{1}{e^{1/x_{\varepsilon }}+1} \right],  \label{h-lim}
\end{equation}
where 
\[
x_{\varepsilon }=x((1+\varepsilon /2)  e^{3i\varepsilon }).
\]
A study of the curves defined by$\Re(\frac{1}{e^{1/x}+1})=0$, that is, 
$e^{1/x}+1=si,s\in \R$   leads to 
\[
x=\frac{1}  {\ln \sqrt{1+s^{2}}+i\arg (-1+si)}  .
\]

\noindent
These curves represent points in the $x$-plane where 
$e^{1/x}+1=$ purely imaginary numbers.We have branches of the curve which is
a consequence of the multi-valueness of $\arg (-1+si)$.  These curves
cluster at $x=0$. \ The point $x_{\varepsilon }=
x((1+\varepsilon/2)e^{3i\varepsilon })$ is in a small vicinity of $x=-i/\pi$.  Note that 
\begin{equation*}
\left.  \frac{dx}{ds} \right|_{x=-i/\pi }=\frac{-i}{\pi ^{2}}.
\end{equation*}
\noindent
This means that the curve has a vertical tangent at $x=-i/\pi$.  In the 
$\zeta$-plane 
\[
\left. \frac{d}{d\varepsilon }(1+\varepsilon /2)
e^{3i\varepsilon }\right| _{\varepsilon =0}=\frac{1}{2}+3i.
\]
 It is easy
to see that the pull back of $x_{\varepsilon }$ does not lie on any of the
branches of the curve 
\[
\Re \left(\frac{1}{e^{1/x}+1} \right)=0
\]
Moreover, since $e^{\pi i}+1=0$, it is not hard to see that 
\[
\left| \Re \left(\frac{1}{e^{1/x_{\varepsilon }}+1} \right)\right| 
\geq   \frac{K}{\varepsilon },
\]
for some positive constant $K$.  This implies from (\ref{h-lim}), 
\[
\left| \widetilde{h}_{n}(1+\varepsilon /2)\right| \geq \frac{K}{\varepsilon }.
\]
Next we estimate 
\begin{equation*}
\max_{  \left| \zeta -(1+\varepsilon /2)\right| =2 \varepsilon }
\left| 
\widetilde{h}_{n}(\zeta )\right| =\max_{\left| \zeta -
(1+\varepsilon/2)\right| =2\varepsilon }\left| \frac{E_{n}(x(\zeta ))\sqrt{n}}
{n!(ex(\zeta))^{n}}\right|,
\end{equation*}
we can still use  (\ref{h-esti})  for the purpose.  This is so 
because the pullback $z(\zeta )$ is still in $\Re(z)<1$. Proceeding similarly and using
 (\ref{h-esti}) ,   we get 
\begin{equation*}
\max_{\left| \zeta -(1+\varepsilon /2)\right| =2\varepsilon }
\left| 
\widetilde{h}_{n}(\zeta )\right| \leq K
\left[(1+O(\varepsilon ))^{n}\sqrt{n}
+O_{\varepsilon }(n^{1-3\alpha }) \right],
\end{equation*}
where the big oh constant in $O(\varepsilon )$ is an absolute constant. 
Further, combining  (\ref{l-ineq})   and  (\ref{l-ineq1}),  we get 
\begin{equation*}
l\leq \frac{K}{\ln 2}\ln \left [(1+O(\varepsilon ))^{n}\sqrt{n}+
O_{\varepsilon}(n^{1-3\alpha })\right].
\end{equation*}
Insert the above in  (\ref{l-ineq2}) , we have
\begin{equation*}
\left| \Im(I_{1})\right| \leq \frac{K}{\ln 2}\ln 
\left[(1+O(\varepsilon))^{n}\sqrt{n}+O_{\varepsilon }(n^{1-3\alpha })\right].
\end{equation*}
Similarly we obtain 
\begin{equation*}
\left| \Im(I_{3})\right| \leq \frac{K}{\ln 2}
\ln \left[(1+O(\varepsilon
))^{n}\sqrt{n}+O_{\varepsilon }(n^{1-3\alpha }) \right].
\end{equation*}
Now from  (\ref{n-esti}) 
\begin{equation*}
\frac{1}{n}N_{n}(0,\theta )-\frac{1}{n}N_{n}(3\varepsilon ,\theta
+\varepsilon _{n}^{\prime })=\frac{1}{n}O(r_{\varepsilon }^{(n)}).
\end{equation*}
Inserting the estimate from $\Im(I_{i}),1\leq i\leq 4$, in the above
and taking the limit sup, we get 
\[
\limsup_{n\rightarrow \infty }\frac{N_{n}(0,\theta )}{n}
\leq \frac{\theta }{2\pi }+O(\varepsilon ).
\]
An analogous lower bound can be obtained, that is, 
\[
\liminf_{n\rightarrow \infty }
\frac{N_{n}(0,\theta )}{n}\geq \frac{\theta }{2\pi }-O(\varepsilon ).
\]
Since $\varepsilon $ can be made arbitrarily small, we conclude 
\[
\lim_{n\rightarrow \infty }\frac{1}{n}N_{n}(0,\theta )=\frac{\theta }{2\pi }.
\]
This ends the proof of Theorem \ref{den-arc}.
\end{proof}

\begin{figure}[ht]
\centerline{\psfig{figure=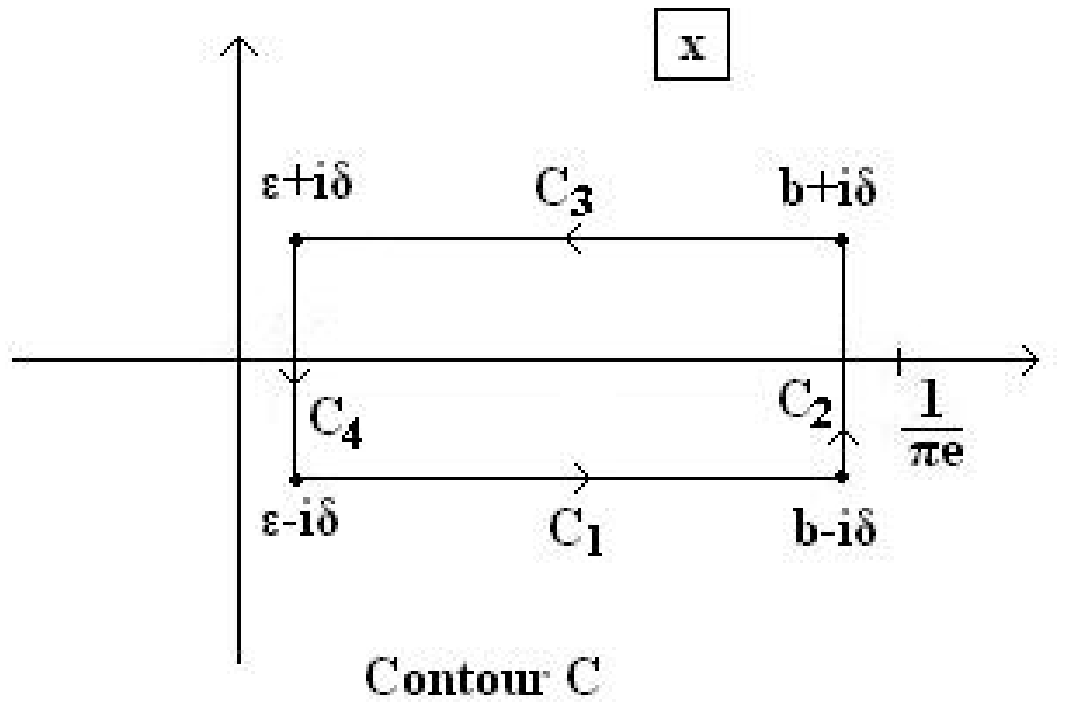, height=3.5in, width=4.0in}}
\caption{}
\label{fig:curveC_6}
\end{figure}

We next calculate the density of zeros in the interval $[0,\frac{1}{\pi e}]$.
 The strategy will be the same to that of Theorem \ref{den-arc}. However,
the technical details are slightly different. We outline the steps below.

\begin{lemma}
(a) For every $\varepsilon >0$, let $N_{\varepsilon }$ be the number of
zeros of $E_{n}(nx)$ in the disc $\{x:\left| x\right| \leq \varepsilon \}$.
\ We show that $N_{\varepsilon }\leq \frac{1}{\ln 2}
(K_{\varepsilon}(1+O(\varepsilon ))^{n})$.

\noindent
(b) For all $0<b<\frac{1}{\pi e}$, we construct the rectangular contour 
$C_{1}\cup C_{2}\cup C_{3}\cup C_{4}=C$ as shown 
in Figure 6.

Let 
$N_{n}(\varepsilon ,b)$ be the number of zeros of $E_{n}(nx)$ enclosed
inside the contour $C$. \ In general, we introduce the notation: let 
$N_{n}(a,b)$ be the number of zeros of $E_{n}(nx)$ enclosed in the
rectangular contour with vertices $a-i\delta ,b-i\delta ,
b+i\delta,a+i\delta $. \ Then for all $\varepsilon >0$ and for all sufficiently small 
$\delta >0$ we have 
\[
\frac{1}{n}N_{n}(\varepsilon ,b)\leq b+O(\delta )+O(\varepsilon )+
O(\tan^{-1}\frac{\delta }{\varepsilon })+
\frac{1}{n\ln 2}\ln (K_{\varepsilon } \left(1+O(\varepsilon )\right)^{n}).
\]
\end{lemma}

 A consequence of the above two lemmas is: $\lim_{n\rightarrow \infty }
\frac{1}{n}N(0,b)=b$. Thus the real roots of $E_{n}(nx)$ are uniformly
distributed in $[-\frac{1}{\pi e},\frac{1}{\pi e}]$.

\begin{proof}
Let us show part (a) first. Since $\frac{E_{n}(nx)}{n!}$ is analytic on the
disc $\{x:\left| x\right| \leq 2\varepsilon \}$, we apply Theorem  \ref{jens}
to get 
\begin{equation}
\left(\frac{2\varepsilon }{\varepsilon }\right)^{N_{\varepsilon }}
=2^{N_{\varepsilon}}\leq 
\frac{\max_{\left| x\right| =2\varepsilon }\left| \frac{E_{n}(nx)}{n!}
\right| }{\left| \frac{E_{n}(0)}{n!}\right| }.  \label{n-esti1}
\end{equation}
$\left| \frac{E_{n}(0)}{n!}\right| $ can be obtained from the integral
representation: 
\begin{equation*}
\frac{E_{n}(0)}{n!}=\frac{1}{\pi i}\oint_{\left| \xi \right| =1}
\frac{1}{(e^{\xi }+1)\xi ^{n+1}}  \,  d\xi .
\end{equation*}
Here caution must be exercised because for even $n$ except $n=0$, 
$E_{n}(0)=0$. This means that $E_{n}(x)$ has a zero at $x=0$ for all even $n\geq 2$.
 Therefore, we must modify the function if we still want to use Jensen's
inequality.  We consider $\frac{\frac{1}{x}E_{n}(nx)}{n!}$ when $n$ is
even.  In this case, an equivalent  Jensen's inequality is: 
\begin{equation*}
2^{N_{\varepsilon }}\leq \frac{\max_{\left| x\right| =2\varepsilon }
\left|  \frac{E_{n}(nx)}{xn!}\right| }    {\left| \lim_{x\rightarrow 0}
\frac{  \frac{1}{x}
E_{n}(nx)     }   {n!}\right| }
\end{equation*}
We will show the former inequality has the main features in the proof, while
the latter can be treated similarly. From integral representation ($n$ is
odd now) we apply the Darboux method \cite{olver}.  In this case the nearest
singularities are $\pm \pi i$.  One shows that 
\begin{eqnarray}
\frac{E_{n}(0)}{n!}
&=&  \frac{1}{\pi i}\oint_{\left| \xi \right| =1}
\left(
\frac{(-1)}{\xi -\pi i}+\frac{(-1)}{\xi +\pi i} \right)\frac{1}{\xi ^{n+1}} \, d\xi +
o(\frac{1}{\pi^{n}})  \nonumber \\
&=&
\frac{2}{(\pi i)^{n+1}}+\frac{2}{(-\pi i)^{n+1}}+o(\frac{1}{\pi ^{n}})
\label{e-esti}
\end{eqnarray}
The corresponding quantity when $n$ is even is 
\begin{equation*}
\frac{\lim_{x\rightarrow 0}\frac{1}{x}E_{n}(nx)}{n!}=\frac{1}{n!}
\lim_{x\rightarrow 0}\frac{1}{x}
\left(
\frac{2}{2\pi i}
\oint_{\left| \xi \right|=1}
\frac{e^{nx\xi }}{(e^{\xi }+1)\xi ^{n+1}}  \,  d\xi 
 \right)
\end{equation*}
Now we know that $\oint_{\left| \xi \right| =1}
\frac{1}{(e^{\xi }+1)\xi^{n+1}} \, d\xi =0$ ($n$ is even here).  
We use the above to rewrite 
$\lim_{x\rightarrow 0}\frac{1}{x}E_{n}(nx)$. Thus 
\begin{eqnarray*}
\lim_{x\rightarrow 0}\frac{1}{x}E_{n}(nx)
&=& \frac{1}{n!}\lim_{x\rightarrow 0}
\frac{2}{2\pi i}
\oint_{\left| \xi \right| =1}
\frac{\frac{1}{x}(e^{nx\xi }-1)}{(e^{\xi }+1)\xi ^{n+1}}
\, d\xi  \\
&=&
\frac{1}{n!}\frac{2}{2\pi i}\oint_{\left| \xi \right| =1}
\frac{n\xi }{(e^{\xi }+1)\xi ^{n+1}}  \,  d\xi  \\
&=&
\frac{1}{(n-1)!}\frac{2}{2\pi i}\oint_{\left| \xi \right| =1}
\frac{1}{(e^{\xi }+1)\xi ^{n}} \, d\xi \quad (n\text{ even})
\end{eqnarray*}
which is still of the same feature as in  (\ref{e-esti}). In this way we see it
does not matter whether $n$ is odd or even, they can be handled in a similar
way. Now we consider 
\begin{equation*}
\max_{\left| x\right| =2\varepsilon }\left| \frac{E_{n}(nx)}{n!}\right| .
\end{equation*}
We use Proposition \ref{m-kprop} with $\mu _{1}$ chosen sufficiently large so
that $\frac{1}{(2\mu _{1}+1)\pi }<2\varepsilon <\frac{1}{(2\mu _{1}-1)\pi }$.
Here the additional assumption that the arbitrarily small number 
$2\varepsilon $ is not of the form $\frac{1}{(2_{\mu _{1}}+1) \pi },  m\geq 0$,
will not hurt the arguments we present here.
Then with $\mu =\mu _{1}-1$ in
Proposition \ref{m-kprop} we get 
\begin{equation*}
\frac{E_{n}(nx)}{n!}=M_{n,\mu _{1}-1}(x)+K_{n,\mu _{1}-1}(x).
\end{equation*}
The asymptotics for $M_{n,\mu _{1}-1}(x)$ is 
\[
M_{n,\mu _{1}-1}(x)
=\sqrt{\frac{2}{\pi }}
\left(  (xe)^{n}F_{\mu _{1}-1}(\frac{1}{x})
\frac{1}{\sqrt{n}x}       \right) \, 
\left(1+O(\frac{1}{n}) \right)
\]
and the asymptotics for $K_{n,\mu _{1}-1}(x)$ comes from applying 
Proposition \ref{szego2}.
This is so because $\left| x(2k+1)\pi i\right| \leq \left| x\right|
(2\mu_{1}-1)\pi $ for all $0\leq k\leq \mu _{1}-1$. \ When $\left| x\right|
=2\varepsilon $, we find $\left| x(2k+1)\pi i\right| \leq 2\varepsilon 
(2\mu_{1}-1)\pi <1$. This means the Proposition \ref{szego2} is applicable.
Furthermore, the largest term comes from $k=0$.  Thus 
\begin{eqnarray}
\lefteqn
{
\frac{E_{n}(nx)  }{n!}
=(xe)^{n}
 [
\sqrt{\frac{2}{\pi }}
F_{\mu _{1}-1}(\frac{1}{x} )  
 \frac{1}{\sqrt{n}x}    \left(1+O(\frac{1}{n}) \right) 
}   \nonumber  \\
&& \qquad + \,
\frac{2g^{n}(x)}{\pi i} \,  (1+o(1))+\frac{2g^{n}(-x)}{(-\pi i)} \, (1+o(1)) 
 ]
\label{m-kequ1}
\end{eqnarray}
Note that on $\left| x\right| =2\varepsilon$, $x$  is in the interior of the
rotated Szeg\"{o} curve $\left| x\pi ie^{1-x\pi i}\right| =1$.  Hence 
$\left| g(x)\right| >1$, $\left| g(-x)\right| >1$.
On the upper semicircle of 
$\left| x\right| =2\varepsilon$, $\left| g(-x)\right| >\left| g(x)\right|$, 
 while on
the lower semicircle of $\left| x\right| =2\varepsilon$, 
$\left| g(x)\right|>\left| g(-x)\right| $, and on the real axis 
$\left| g(x)\right| =\left|g(-x)\right|$.  No matter what we always have 
\begin{equation*}
\left| g(x)\right| \leq \frac{\left| e^{x\pi i}\right| }
{e\left| x\right| \pi }=\frac{1+O(\varepsilon )}{e\left| x\right| \pi },\text{for }
\left|  x\right| =2\varepsilon.
\end{equation*}
Similarly $\left| g(-x)\right| \leq \frac{1+O(\varepsilon )}
{e\left|x\right| \pi }$.  Hence from  (\ref{m-kequ1})  we conclude 
\[
\left| \frac{E_{n}(nx)}{n!}\right| \leq 
\frac{K_{\varepsilon}(1+O(\varepsilon ))^{n}}{\pi ^{n}}
\]
for some absolute constant $K_{\varepsilon }$ depending only on $\varepsilon$.
Using (\ref{n-esti1}) and (\ref{e-esti}), we get 
\begin{eqnarray*}
2^{N_{\varepsilon }}
&\leq&
 \frac{\max_{\left| x\right| =2\varepsilon }
\left|  \frac{E_{n}(nx)}{n!}\right| }{\left| \frac{E_{n}(0)}{n!}\right| } \\
&\leq &
\frac{\frac{K_{\varepsilon }(1+O(\varepsilon ))^{n}}{\pi ^{n}}}{(\frac{A}{\pi^{n}})} \\
&=&
K_{\varepsilon }(1+O(\varepsilon ))^{n}
\end{eqnarray*}
Taking logarithms we get 
\begin{equation*}
N_{\varepsilon }  
\leq \frac{1}{  \ln 2}
\ln \left(  K_{\varepsilon }
(1+O(\varepsilon)   \right)^{n}.
\end{equation*}
This ends the proof of part (a).

We now prove part (b). 
\begin{equation*}
N(\varepsilon ,b)=\frac{1}{2\pi }
\Im
\left(
\oint_{C} \frac{h_{n}^{\prime }(x)}   {h_{n}(x)}  \,  dx 
\right)
\end{equation*}
where $h_{n}(x)=\frac{E_{n}(nx)\sqrt{n}}{n!(ex)^{n}}$. 
Decomposing $C$ into 
$C_{1},  C_{2},  C_{3}$,  and $C_{4}$ as shown in Figure 6, we get 
\begin{equation}
N(\varepsilon ,b)=\frac{1}{2\pi }\Im(J_{1})+\frac{1}{2\pi }
\Im(J_{2})+
\frac{1}{2\pi }\Im(J_{3})+\frac{1}{2\pi }\Im(J_{4}),
\label{n-inte}
\end{equation}

\noindent
where 
\begin{equation*}
J_{i}=\int_{C_{i}}\frac{h_{n}^{\prime }(x)}{h_{n}(x)} \,
dx, \quad 1  \leq i  \leq 4.
\end{equation*}
Let us focus on $J_{1}$ first. \ We infer that (\ref{m-kequ1}) still works for
all $x$ on $C_{1}$ provided $\Im(x)=\delta $ is sufficiently small. 
This is so, for $\left| g(x)\right| =\frac{\left| e^{x\pi i}\right| }
{ e\left| x\right| \pi }=\frac{e^{x\left| x\right| \sin \theta }}
{e\left| x\right| \pi }$, where for $x\in C_{1},x=\left| x\right| 
e^{-i\theta},0<\theta <\pi $. This observation leads to 
\begin{equation}
\frac{E_{n}(nx)\sqrt{n}}{n!(ex)^{n}}
=
 \sqrt{\frac{2}{\pi }}F_{\mu _{1}-1}
( \frac{1}{x})\frac{1}{x} \, \left(1+O(\frac{1}{n}) \right)  + 
\frac{2g^{n}(x)}{\pi i} \, \left(1+o(1) \right),  \label{e-esti1}
\end{equation}

\noindent
that is, $g^{n}(x)$ is the dominant term among the terms in 
$K_{n,\mu_{1}-1}(x)$. Hence $h_{n}(x)\cdot g^{-n}(x)\rightarrow \frac{2}{\pi i}$
uniformly. This implies 
\begin{equation*}
\frac{h_{n}^{\prime }(x)}{h_{n}(x)}-n(\pi i-\frac{1}{x})\rightarrow 0,
\text{  uniformly}
\end{equation*}
Integrating the above along $C_{1}$ we get 
\begin{eqnarray*}
\frac{1}{2\pi n}\Im(J_{1})
&\rightarrow&
\Im \left(\int_{C_{1}}(\pi i-  \frac{1}{x})  \, dx \right) \\
&=&
\frac{1}{2\pi }\Im
\left[\pi i(b-i\delta -\varepsilon +i\delta )-
\ln  \frac{b-i\delta }{\varepsilon -i\delta } \right]  \\
&=&
\frac{1}{2\pi }
\left(\pi (b-\varepsilon )-\arg \frac{b-i\delta }
{\varepsilon  -i\delta } \right)  \\
&=&
\frac{b-\varepsilon }{2}+O(\delta )+
O(\tan ^{-1}\frac{\delta }{\varepsilon })
\end{eqnarray*}
that is, 
\begin{equation*}
\lim_{n\rightarrow \infty }\frac{1}{2\pi n}\Im(J_{1})=\frac{b}{2}
+O(\delta )+O(\varepsilon )+O(\tan ^{-1}\frac{\delta }{\varepsilon })
\end{equation*}
Similarly, we obtain 
\begin{equation*}
\lim_{n\rightarrow \infty }\frac{1}{2\pi n}I(J_{3})=\frac{b}{2}+
O(\delta)+O(\varepsilon )+O(\tan ^{-1}\frac{\delta }{\varepsilon })
\end{equation*}
(In this case, $g^{n}(-x)$ becomes dominant, rather than $g^{n}(x)$.)  An
estimate for $\frac{1}{2\pi }\Im(J_{2})$ comes from the observation
that the change of arguments for $g(x)$ on the vertical segment 
$\overline{b-i\delta ,b}$ is of order $O(\delta )$  (\ref{e-esti1}). Hence, 
\begin{equation*}
\frac{1}{2\pi n}\Im(J_{2})=O(\delta ).
\end{equation*}
Similar to (\ref{n-esti1})  we can likewise prove that 
\begin{equation*}
\frac{1}{2\pi n}\Im(J_{4})\leq \frac{1}{n\ln 2}
\ln (K_{\varepsilon}(1+O(\varepsilon ))^{n}).
\end{equation*}
Inserting all these estimates into (\ref{n-inte})  we get 
\[
\frac{1}{n}N_{n}(\varepsilon ,b)
\leq
 b+O(\delta )+O(\varepsilon )+
O(\tan^{-1}\frac{\delta }{\varepsilon })+ 
\frac{1}{n\ln 2}\ln (K_{\varepsilon }  \,  (1+O(\varepsilon ))^{n}).
\]
Again taking (\ref{n-esti1})  into account we get 
\[
\frac{1}{n}N_{n}(0,b)
\leq   b+O(\delta )+O(\varepsilon )+
O(\tan ^{-1}\frac{ \delta }{\varepsilon })+ 
\frac{2}{n\ln 2}\ln (K_{\varepsilon }(1+O(\varepsilon ))^{n}).
\]
This implies 
\[
\limsup_{n\rightarrow \infty }  
\frac{1}{n}N_{n}(0,b)\leq b+
O(\delta)+O(\varepsilon )+  
O(\tan ^{-1}\frac{\delta }{\varepsilon })+\ln \left(1+O(\varepsilon ) \right).
\]
We get a similar lower bound for 
\[
\lim \inf_{n\rightarrow \infty } 
\frac{1}{n}N_{n}(0,b)\geq b-O(\delta)-O(\varepsilon )-
O(\tan ^{-1}\frac{\delta }{\varepsilon })-\ln \left(1+O(\varepsilon ) \right).
\]
But the above is true for all $\varepsilon >0$ and all $\delta >0$, hence 
\begin{equation*}
\lim_{n\rightarrow \infty }\frac{1}{n}N_{n}(0,b)=b.
\end{equation*}
In this way we have established all claims we have made in this paper.
\end{proof}

\end{document}